\documentclass[sn-mathphys,Numbered]{sn-jnl}

\usepackage{graphicx}%
\usepackage{multirow}%
\usepackage{amsmath,amssymb,amsfonts}%
\usepackage{amsthm}%
\usepackage{mathrsfs}%
\usepackage[title]{appendix}%
\usepackage{xcolor}%
\usepackage{textcomp}%
\usepackage{manyfoot}%
\usepackage{booktabs}%
\usepackage{algorithm}%
\usepackage{algorithmicx}%
\usepackage{algpseudocode}%
\usepackage{caption}
\usepackage{listings}%
\usepackage{mathtools}
\usepackage{doi}
\usepackage{lmodern}

\usepackage{float}
\restylefloat{table}
\usepackage{epstopdf}
\usepackage{subfigure}
\usepackage{ifthen}
\numberwithin{equation}{section}
\usepackage{helvet}
\usepackage{fourier}
\usepackage{array}
 \usepackage{makecell}

\usepackage{graphics}
\usepackage{pstricks}
\usepackage{comment}
\definecolor{codegreen}{rgb}{0,0.6,0}
\definecolor{codegray}{rgb}{0.5,0.5,0.5}
\definecolor{codepurple}{rgb}{0.58,0,0.82}
\definecolor{backcolour}{rgb}{0.95,0.95,0.92}
\lstdefinestyle{matlab}{
    backgroundcolor=\color{backcolour},
    commentstyle=\color{codegreen},
    keywordstyle=\color{magenta},
    numberstyle=\tiny\color{codegray},
    stringstyle=\color{codepurple},
    basicstyle=\ttfamily\footnotesize,
    breakatwhitespace=false,
    breaklines=true,
    captionpos=b,
    keepspaces=true,
    numbers=left,
    numbersep=5pt,
    showspaces=false,
    showstringspaces=false,
    showtabs=false,
    tabsize=2
}

\theoremstyle{definition}

\newtheorem{defn}{Definition}[section]

\DeclareMathOperator{\card}{card}
\DeclareMathOperator{\conv}{conv}
\DeclareMathOperator{\spn}{span}
\newcommand{\md}{{\rm d}}
\newcommand{\dx}{\, \md x}

\newcommand{\dgamma}{\, \md \gamma}

\theoremstyle{thmstyleone}%
\theoremstyle{thmstyletwo}%
\newtheorem{example}{Example}%
\newtheorem{remark}{Remark}%

\theoremstyle{thmstylethree}%

\begin{document}

\title[Article Title]{Vectorized implementation of primal hybrid FEM in MATLAB}


\author[1]{\fnm{Harish} \sur{Nagula} \sur{Mallesham}}\email{mat21h.nagula@stuiocb.ictmumbai.edu.in}

\author[2]{\fnm{Kamana} \sur{Porwal}}\email{ kamana@maths.iitd.ac.in}

\author[3,4]{\fnm{Jan} \sur{Valdman}}\email{ jan.valdman@utia.cas.cz }

\author*[1]{\fnm{Sanjib Kumar} \sur{Acharya}}\email{sk.acharya@iocb.ictmumbai.edu.in}

\affil[1]{\orgdiv{
Institute of Chemical Technology Mumbai }, \orgname{IndianOil Odisha Campus}, \orgaddress{ \city{Bhubaneswar}, \postcode{751013}, \state{Odisha}, \country{India}}}

\affil[2]{\orgdiv{Department of Mathematics}, \orgname{Indian Institute of Technology Delhi}, \orgaddress{ \city{New Delhi}, \postcode{110016}, \state{Delhi}, \country{India}}}

\affil[3]{\orgdiv{
Faculty of Information Technology}, \orgname{Czech Technical University in Prague}, \orgaddress{Thákurova 9, 16000 Prague}, \country{Czech Republic}}

\affil[4]{\orgdiv{The Czech Academy of Sciences}, \orgname{Institute of Information Theory and Automation}, \orgaddress{ \city{Pod vod\'{a}renskou v\v{e}\v{z}\'{i} 4, 18208, Prague 8}, \country{Czech Republic}}}


\abstract{We present efficient MATLAB implementations of the lowest-order primal hybrid finite element method (FEM) for linear second-order elliptic and parabolic problems with mixed boundary conditions in two spatial dimensions. We employ the Crank-Nicolson finite difference scheme for the complete discrete setup of the parabolic problem. All the codes presented are fully vectorized using matrix-wise array operations. Numerical experiments are conducted to show the performance of the software.
  \\}




\maketitle

\section{Introduction}\label{sec1}
MATLAB is a popular computing platform in academia and industry, with many built-in functions and toolboxes developed by Mathworks Computer Software Corporation. However, codes containing for loops, particularly the assembly of finite elements of stiffness and mass matrices in MATLAB are extremely slow compared to C and FORTRAN (cf. \cite{alberty1999remarks,carstern}). To overcome this issue, an efficient and flexible MATLAB assembly procedure of nodal element stiffness and mass matrices for 2D and 3D was introduced by Rahman and Valdman \cite{rahman2013fast}. They vectorized for-loops by extending the element-wise array operations into matrix-wise array operations, where the array elements are matrices rather than scalars, resulting in a faster and more time-scalable algorithm. This idea is exploited for the MATLAB assembly of 2D and 3D edge elements by Anjam and Valdman in \cite{anjam2015fast}. In \cite{funken2011efficient}, Funken et al. adopted a similar approach while assembling the element stiffness matrices. For more detailed literature, we refer to Cuvelier et al., see \cite{cuvelier2016efficient}, Moskovka and Valdman \cite{moskovka2022fast}  and the references therein.

In this paper, we present an efficient MATLAB implementation procedure of the primal hybrid finite element method, which is based on an extended variational principle introduced by Raviart and Thomas in \cite{raviart} for elliptic problems. Acharya and Porwal \cite{acharya,kamna} extended this technique to parabolic second- and fourth-order problems and established the optimal order error estimates.  The benefit of this method is that two unknown variables (primal and hybrid) can be computed simultaneously without losing accuracy. The constrained space using Lagrange multipliers leads to a non-conforming space; for instance, the lowest-order primal hybrid FEM (using unconstrained space) is equivalent to the Crouzeix-Raviart nonconforming (using constrained space) finite element space. It is also called the generalized nonconforming method.  For a detailed review of the literature and applications, see \cite{acharya,kamna,peter, Aposteriori} and the references therein. To the best of the author's knowledge, no work has been reported on the implementation procedure of this method.

The main contributions are as follows.

\begin{itemize}
\item Edge generation function edge.m (cf. \cite{carstern}) and a uniform mesh refinement function based on red-refinement (cf. \cite{redrefine,bey2000simplicial}) have been vectorized using matrix-wise array operations.
\item We present a vectorized implementation of the lowest-order primal hybrid FEM for a second-order general elliptic problem with mixed boundary conditions.
\item Also, we implement this method for the general second-order parabolic problem with mixed boundary conditions using the Crank-Nicolson finite difference scheme in the temporal direction.

\item 
Numerical experiments have been conducted, and run-time has been presented. The experiments were carried out using MATLAB R2020a on a computer equipped with an $\mathrm{i}5-10210\mathrm{U}$ CPU operating at a frequency of 1.60GHz-2.10 GHz, with 16GB  RAM and a $\mathrm{x}64$-based processor with 1TB of system memory. 
A linear regression analysis shows that the codes have an approximate optimal-order time scale.  

\item The software can be downloaded  through the following link: 
\url{https://www.mathworks.com/matlabcentral/fileexchange/136359}.  

\end{itemize}

The outline of this paper is as follows: In Section \ref{sec:triangulation}, we present the required data structure to present the primal hybrid FEM and its implementation procedures. The primal hybrid finite element algebraic formulation for the elliptic problem and its MATLAB implementation have been discussed in Section \ref{sec:ellipticproblem}. Implementation of the primal hybrid FEM with the Crank-Nicolson scheme in the time direction for parabolic problems is presented in Section \ref{sec:parabolicproblem}. 
We give our concluding remarks in Section \ref{conclusion}.

\section{Triangulation and geometric structure}\label{sec:triangulation}

Let $\Omega=(0,1)^2$ be the computational domain with boundary $\Gamma=\overline{\Gamma}_{D}\cup\overline{\Gamma}_N$, where $\overline{\Gamma}_D$  and $\Gamma_N=\Gamma/\overline{\Gamma}_{D}$ are the Dirichlet and Neumann boundary, respectively. We use the Sobolev spaces $H^{m}(\Omega)$, $H(\text{div},\Omega)$ in the domain $\Omega$, $m\in\mathbb{R}_{\pm}$ \cite{ciar}.
With the mesh parameter $h$, let $\mathcal{T}_{h}$ be a regular family of triangulations of the set $\overline{\Omega}$ in the sense
of Ciarlet \cite{ciar} with shape regular triangles $T$ whose diameters $h_{T}\leq h$ are such that
$
\overline{\Omega}=\bigcup_{T\in\mathcal{T}_{h}}\overline{T}.
$

\begin{figure}[h!]
\centering
\includegraphics[ width=8cm]{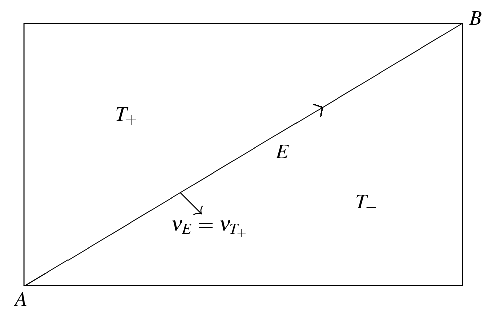}
\centering
\caption{A pair of adjacent triangles $(T_+, T_-)$ with a common edge $E=\partial T_+\cap T_-$, initial node $A$, end node $B$ and outward unit normal $\nu_E=\nu_{T_+}$.}
\label{jump}
\end{figure}

For any triangle $T$, let $\partial T$ denote the boundary of $T$ and let $\nu_T$ be the outward unit normal to $\partial T$.
Let $$\mathcal{E}^{h}=\mathcal{E}^{h}_{\Omega}\cup\mathcal{E}^{h}_{D}\cup\mathcal{E}^{h}_{N}$$ be the set of all edges $E$ of  $\mathcal{T}_{h}$
 where $\mathcal{E}^{h}_{\Omega}$ denotes the set of all interior element edges, 
 $$\mathcal{E}^{h}_{D}=\{E\in\mathcal{E}^{h}:E\subset\overline{\Gamma}_D\}$$ denotes the set of edges on the Dirichlet boundary  and 
 $$\mathcal{E}^{h}_{N}=\{E\in\mathcal{E}^{h}:E\subset\overline{\Gamma}_N\}$$    denotes the set of edges on the Neumann boundary.
The following notations will be used throughout the article:

\noindent\begin{minipage}{.5\linewidth}
\begin{align*}
  \mathcal{N}&= \text{set of vertices of element}~ T \in \mathcal{T}_h,\\
m_E&=\text{mid-point of} ~E \in \mathcal{E}^h,\\
\mathcal{N}_m &= \text{set of mid-points $m_E$ of the edges}~ E \in \mathcal{E}^h,\\
\mathcal{N}_m(\Gamma_D)&=\text{set of midpoints of the edges}~ E \in \mathcal{E}^{h}_{D},\\
\mathcal{N}_m(\Gamma_N)&= \text{set of midpoints of the edges}~ E \in \mathcal{E}^{h}_{N},\\
\card(S)&=\text{cardinality of the set $S$},\\\ \conv \{A,B\}&=\text{Convex hull of} ~A~ \text{and}~B,\\
\end{align*}
\end{minipage}%
\begin{minipage}{.5\linewidth}
\begin{align*}
  \texttt{nE}&=  \card(\mathcal{T}_h), \\
N&= 3 \times \texttt{nE}, \\
L&= \card(\mathcal{E}^h_{\Omega} \cup \mathcal{E}^h_{D}),  \\
\texttt{nrEdges}&= \card( \mathcal{E}^h), \\
N_{D}&= \card( \mathcal{E}^h_{D}), \\
N_{\Omega}&= \card(\mathcal{E}^h_{\Omega}),\\
N_{n}&= \card(\mathcal{E}^h_{N}).\\
\end{align*}
\end{minipage}

 The intersection of two distinct triangles of $\mathcal{T}_h$ is either empty or non-empty. They share one complete edge $E=\conv\{z_1,z_2\}$ or a node $z_1$ in the non-empty case. \figurename{~\ref{jump}} depicts two adjacent triangles $T_+$ and $T_-$ with a common edge $E=\conv\{A,B\}$. Let  $\mathcal{N}= \{z_1,\dots,z_{\mathrm{card}(\mathcal{N})}\}$ be the set of vertices (nodes) $z_l,~1\leq l\leq \card(\mathcal{N})$ with coordinates $(x_l,y_l) \in \mathbb{R}^2$. For an edge $E = \conv \{ z_i,z_j \}$ such that $z_i$ and $z_j$ are ordered counterclockwise, its outward unit normal  defined as: 
 $$\nu_{E}=(y_j-y_i, x_j-x_i)/|E|,$$ where $|E|$ is the length of $E$. We define the jump of a scalar-valued function $v$ on any $E\in\mathcal{E}^{h}$ such that
$E=\partial T_{+}\cap\partial T_{-}$ ($T_{+}$ and $T_{-}$ are adjacent triangles in $\mathcal{T}_{h}$)  as:
$${[\![v]\!]}_{E}= {({v|}_{T_{+}})|}_{E}-{({v|}_{T_{-}})|}_{E}$$
and for a boundary edge $E$ such that
 $E=\partial T\cap\Gamma$, where $T\in\mathcal{T}_h$, we define  ${[\![v]\!]}_{E}={v|}_{T}\nu_{E}$.

\begin{figure}[H]
\centering
\includegraphics[width=8cm]{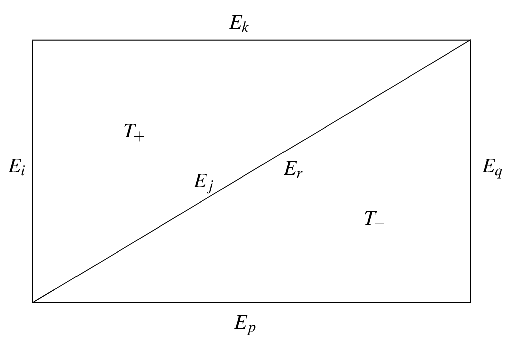}
\centering
\caption{Depiction of a general numbering of edges in a mesh with two adjacent triangles.}
\label{edgenumbering}
\end{figure}

 Let $i,j,k$ represent the global numbering of nodes of a triangle $T=\conv\{z_i,z_j,z_k\}$.
Let  $E_{i}, E_{j}, E_{k}$ be edges  of triangle $T_+$ and $E_{p}, E_{q}, E_{r}$ be edges  of triangle $T_-$ for two adjacent triangular elements $T_+=\conv \{z_i, z_j, z_k\}$ and $T_-=\conv \{z_p, z_q, z_r\}$ in $\mathcal T_h$ as in \figurename{~\ref{edgenumbering}}.

For a triangular element, $T = \conv\{P_1,P_2,P_3\}$, where $$P_1=(x_1,y_1), \quad  P_2=(x_2,y_2), \quad P_3=(x_3,y_3)$$ are vertices of $T$, opposite edges are 
$$E_1 = \conv\{P_2, P_3\}, \quad  E_2 = \conv\{P_1, P_3\}, \quad E_3 = \conv\{P_1, P_2\},$$ mid-point of edges are 
$$mE_1=(P_2+P_3)/2,\quad  mE_2=(P_1+P_3)/2, \quad mE_3=(P_1+P_2)/2,$$ and area of the element is denoted as $|T|$.

\subsection{Data representation of the triangulation}

\begin{figure}[H]
\centering
\includegraphics[ width=6.5cm]{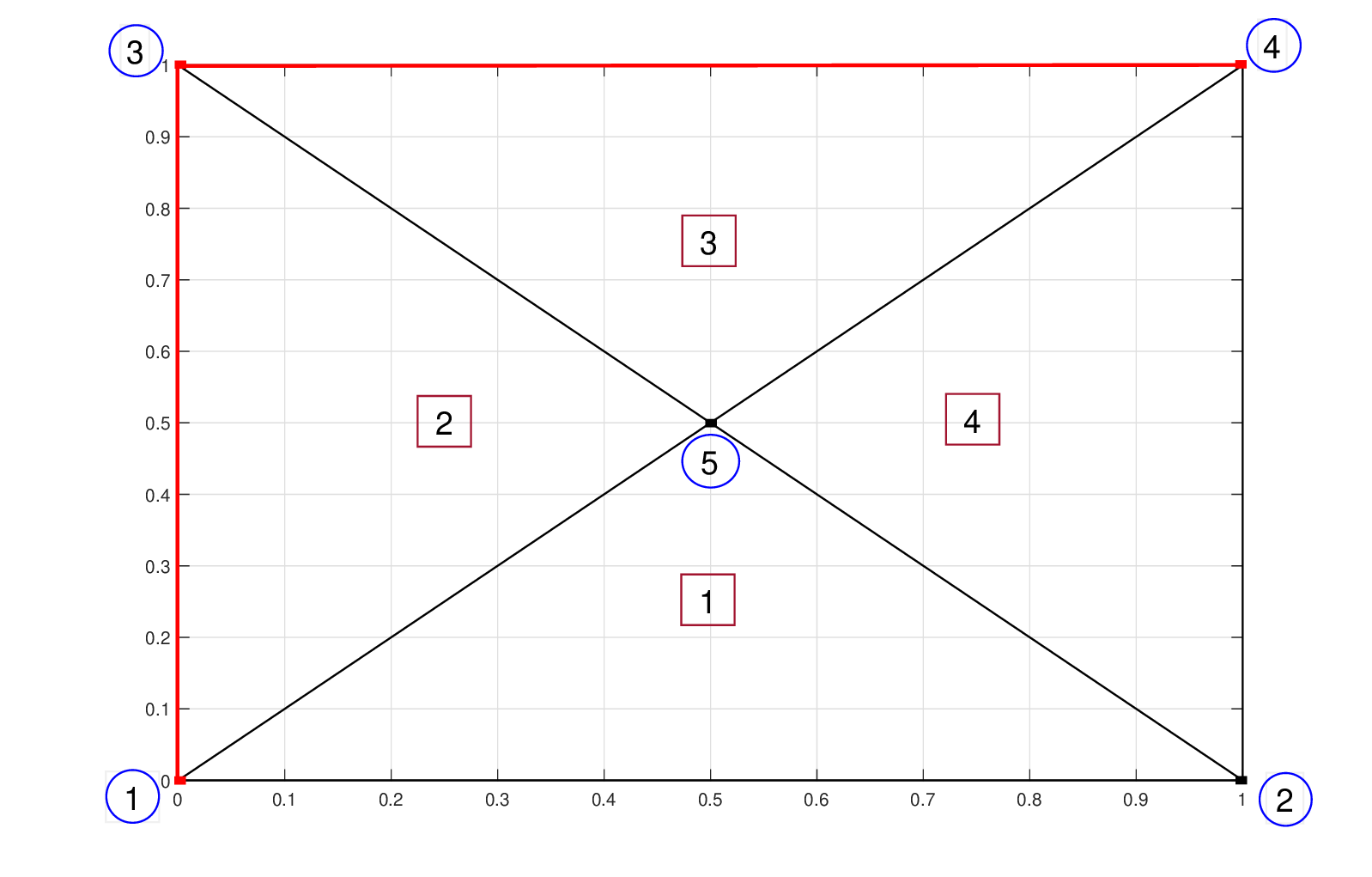}
\includegraphics[width=6.5cm]{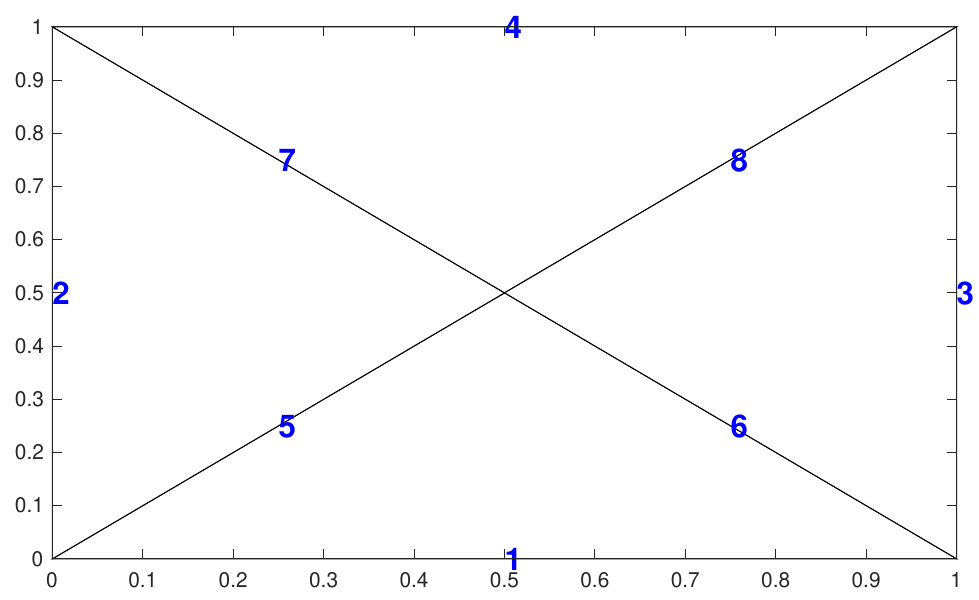}
\centering
\caption{Depiction of initial (level 0) triangulation (left) and the corresponding edges generated (right).
}
\label{InitialTriangulation}
\end{figure}

We present the data structure of the initial triangulation of $\Omega$ with triangles in Table \ref{Table1} and \figurename{~\ref{InitialTriangulation}}.
 The nodes and elements are stored in the rows of $"\texttt{coordinates.dat}"$ and $"\texttt{elements.dat}"$, respectively. The information about the boundary edges in $\mathcal{E}_D^{h}$ and $ \mathcal{E}_{N}^{h}$ are stored in  $"\texttt{Dirichlet.dat}"$ and $"\texttt{Neumann.dat}"$, respectively.
\begin{table}[ht]
\footnotesize
\captionsetup{font=footnotesize}
            \begin{minipage}[b]{0.15\linewidth}
           \centering
                \begin{tabular}[t]{c c}
                 \multicolumn{2}{c}{\texttt{coordinates.dat}}
                 \vspace{0.2cm} \\
                    \hline
                    \rule{0pt}{1.5\normalbaselineskip}
                    \vspace{0.2cm}
                    0&   0\\
                    \vspace{0.2cm}
                    1&0 \\
                    \vspace{0.2cm}
                    0&1\\
                    \vspace{0.2cm}
                    1&1 \\
                    \vspace{0.2cm}
                    0.5&0.5\\  \hline 
                \end{tabular}
                \end{minipage}
                \hspace{1.5cm}
                \begin{minipage}[b]{0.15\linewidth}
    \centering
                \begin{tabular}[t]{c c c}
                \multicolumn{3}{c}{\texttt{elements.dat}} \vspace{0.2cm} \\
                    \hline
                    \rule{0pt}{1.5\normalbaselineskip}
                    \vspace{0.2cm}
                    5&1&2\\
                    \vspace{0.2cm}
                    5&3&1 \\
                    \vspace{0.2cm}
                    5&4&3\\
                    \vspace{0.2cm}
                    5&2&4\\  \hline
                \end{tabular}
                \end{minipage}
                \hspace{1cm}
                \begin{minipage}[b]{0.15\linewidth}
    \centering
                \begin{tabular}[t]{c c}
                    \multicolumn{2}{c}{\texttt{Dirichlet.dat}} \vspace{0.2cm}\\
                    \hline
                    \rule{0pt}{1.5\normalbaselineskip}
                    \vspace{0.2cm}
                    1&2\\
                    \vspace{0.2cm}
                    2&4 \\    \hline
                \end{tabular}
                \end{minipage}
                \hspace{1cm}
                \begin{minipage}[b]{0.05\linewidth}
    \centering
                \begin{tabular}[t]{c c}
                    \multicolumn{2}{c}{\texttt{Neumann.dat}} \vspace{0.2cm}\\
                    \hline
                    \rule{0pt}{1.5\normalbaselineskip}
                    \vspace{0.2cm}
                    4&3\\
                    \vspace{0.2cm}
                    3&1 \\ \hline
                \end{tabular}
                \end{minipage}
                \caption{\texttt{coordinates.dat}, \texttt{elements.dat}, \texttt{Neumann.dat} (consists of 2 boundary edges which are plotted in red) and  \texttt{Dirichlet.dat} (remaining 2 boundary edges) data for the triangulation displayed in \figurename{ \ref{InitialTriangulation}.}}
                \label{Table1}
\end{table}
\subsection{Vectorization of $\mathbf{edge.m}$}
In this subsection, we present a modified version of the function $\mathbf{edge.m}$ (whose inputs are \texttt{elements} and \texttt{coordinates} and outputs are $\texttt{nodes2element}(\texttt{n2el})$, $\texttt{nodes2edge}(\texttt{n2ed})$ and $\texttt{edge2element}(\texttt{ed2el})$ (cf. Bahriawati and Carstensen, \cite{carstern}). The following example intuitively explains the meaning of output matrices. \\

\begin{example}
Given the triangulation of Fig.  \ref{InitialTriangulation}, we have the following full versions of sparse matrices
\begin{equation}
\texttt{nodes2element} =
\begin{pmatrix}
     0    & 1  &   0    & 0  &   2 \\
     0   &  0   &  0  &   4  &   1 \\
     2  &   0  &   0   &  0   &  3 \\
     0  &   0  &   3   &  0  &   4 \\
     1  &   4  &   2  &   3  &   0    
\end{pmatrix}, 
\qquad 
\texttt{nodes2edge}=
\begin{pmatrix}
   0  &   1  &   2  &   0  &   5 \\
     1   &  0  &   0  &   3  &   6 \\
     2   &  0  &   0  &   4  &   7 \\
     0  &   3  &   4  &   0   &  8 \\
     5  &   6  &   7   &  8  &   0
\end{pmatrix}.
\end{equation}
The 1st row of \texttt{nodes2element} tells us that
the directed edge $1 \rightarrow 2$ belongs to the element 1 and
the directed edge $1 \rightarrow 5$ belongs to the element 2. 
Note the directed edge $3 \rightarrow1$ is stored in the 3rd raw and belongs to the element 2. \\

The matrix \texttt{nodes2edge} provides actual numbers of edges. The first raw reveals that the undirected edge $1 \leftrightarrow 2$ has number 1, the undirected edge $1 \leftrightarrow 3$ has number 2, and the undirected edge $1 \leftrightarrow 5$ has number 5. \\

Finally, the matrix
$\texttt{edge2element}$ contains nodes belonging to edges (the first two columns) and the connection of edges to elements (the remaining two columns). 
\begin{equation}
\texttt{edge2element}=
\begin{pmatrix}
     1  &   2  &   1   &  0  \\
     3   &  1  &   2  &   0   \\
     2   &  4  &   4  &   0 \\
     4  &   3  &   3   &  0   \\
      5  &   1  &   1  &   2    \\
     2   &  5  &   1  &   4  \\
     5  &   3  &   2   &  3   \\   
     5   &  4   &  3   &  4   \\    
\end{pmatrix}.
\end{equation}
The fifth row tells us that the edge $5 \leftrightarrow 1$ belongs to the elements 1, 2. Clearly, the rows of $\texttt{edge2element}$ with the fourth column entry equal to zero describe boundary edges that belong to one element only. 
\end{example}
\medskip

By removing for-loops of  $\mathbf{edge.m}$ function, we provide a modified function: \\ \\
\noindent
[$\texttt{nodes2element, nodes2edge, edge2element}$]=$\mathbf{edge\_vec}(\texttt{elements, coordinates})$ \\ \\
\noindent
in which the interior edges (\texttt{intE}), exterior edges (\texttt{extE}), local edge numbering data for $T_{+}$ elements (\texttt{led}), and local edge numbering data for $T_{-}$ elements (\texttt{ledTN}) can additionally be output.

\subsubsection{\texttt{nodes2element}}
The $\texttt{nodes2element}$ function describes an element's global number as a function of its two vertices,  with dimension $\card(\mathcal{N})\times\card(\mathcal{N})$ defined as,
\begin{equation}\label{eq:n2el}
        \texttt{nodes2element}~(k,l) = \begin{cases} j & \mbox{if } (k,l) ~\mbox{are numbers of the nodes of  element number} ~ j;\\
 0 & \mbox{otherwise.} \end{cases}
\end{equation}

\begin{figure}[H]
\centering
\includegraphics[width=10.5cm]{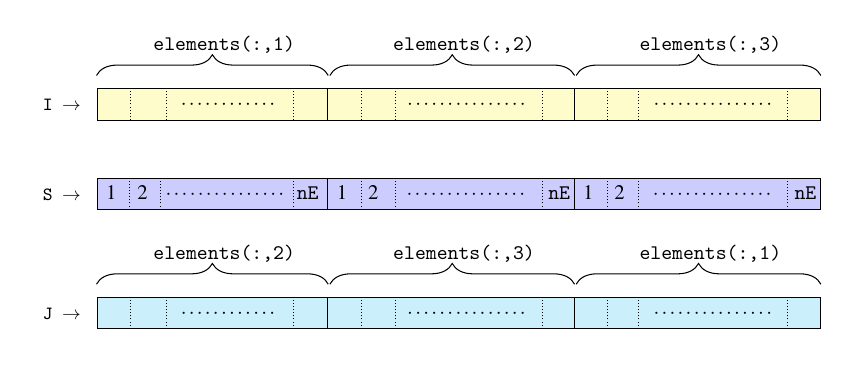}
\centering
\caption{Geometrical representation of \texttt{I}, \texttt{J}, and \texttt{S} vectors in MATLAB.}
\label{fig:node2element_diagram}
\end{figure}

For the matrix \texttt{nodes2element}, the vectorized computation in MATLAB is as follows:
\lstset{language=Matlab,style=matlab}
\begin{lstlisting}[language=Matlab, caption={Vectorization of the \texttt{nodes2element} matrix.},label=lst:vec_n2el]
nE=size(elements,1); %number of elements
nC=size(coordinates,1); % number of nodes
I=elements(:);
J=reshape(elements(:,[2 3 1]), numel(elements), 1);
S=reshape([1:nE;1:nE;1:nE]', numel(elements), 1);
nodes2element=sparse(I, J, S, nC, nC);
\end{lstlisting}
Here \texttt{nodes2element} is a sparse matrix from the triplets \texttt{I}, \texttt{J}, and \texttt{S} such that, $
\texttt{nodes2element}(~\texttt{I}(k),~\texttt{J}(k))= \texttt{S}(k)
$, where \texttt{I}, \texttt{J}  and \texttt{S} are column-vectors as shown in \figurename{~\ref{fig:node2element_diagram}}.
\subsubsection{\texttt{nodes2edge}}
\vspace{-0.1cm}
The  sparse matrix $\texttt{node2edge}$ describes the global number of edges, with dimension $\card(\mathcal{N})\times\card(\mathcal{N})$ defined as,
\begin{align}\label{eq:n2ed}
        \texttt{nodes2edge}~(k,l) = \begin{cases} j ~& \mbox{if } j-\mbox{th} \mbox{  edge is connected between the}\\
        & \mbox{two nodes with global numbering} ~ (k,l);  \\
 0 ~& \mbox{otherwise.} \end{cases}
\end{align}
The MATLAB code for the calculation of the matrix \texttt{nodes2edge} is as follows:
\begin{lstlisting}[language=Matlab, caption={Vectorization of the \texttt{nodes2edge} matrix.},label=lst:vec_n2ed]
B=nodes2element+nodes2element';
[i,j]=find(triu(B));
nodes2edge=sparse(i,j,1:size(i,1),nC,nC);
nodes2edge=nodes2edge+nodes2edge';
\end{lstlisting}
In the above code, the variables \texttt{i} and \texttt{j} represent the row and column indices for the nonzero entries of the upper triangular part of the matrix \texttt{B}. The number of edges is denoted by \texttt{nrEdges=size(i,1)}.

\subsubsection{$\texttt{edge2element}$}
 The row number $p1$ of the matrix \texttt{edge2element} represents the edge's initial node $k$ and end node $l$, as well as the numbers $m, n$ of elements $T_+ , T_-$ that share the edge, and the matrix \texttt{edge2element} has a dimension of \texttt{nrEdges} $\times 4$.
\begin{equation*}
        \texttt{edge2element}~(p1,[1,2])=[k ~~ l] ~~~~~~(k,l)  ~\mbox{are initial node and end node of edge}~p1.
\end{equation*}
The \figurename{~\ref{jump}} convention is used to specify the neighbor elements $T_+, T_-$ with indices $m$ and $n$. The entry of the matrix $\texttt{edge2element}(p1,3)$ defines the element $T_+$, which also specifies the orientation of the unit normal $\nu_E=\nu_{T_+}$ of the edge $E = T_+ \cap T_-$. 
If the edge $E \in \mathcal{E}^h_{D} \cup \mathcal{E}^h_{N}$, then the fourth entry $\texttt{edge2element}(p1,4)$) is zero, indicating that $\nu_E$ lies on $\partial\Omega$.
\begin{equation}\label{eq:ed2el}
 \texttt{edge2element}~(p1,[3,4]) = \begin{cases} [m~n]  ~& \mbox{if the common edge } p1 \mbox{ belongs to element} ~ m ~\mbox{and}~ n; \\
 [m~0]  ~& \mbox{if the boundary edge}~ p1 \mbox{ belongs to element}~m .\end{cases}
\end{equation}
The MATLAB code for the vectorized calculation of the matrix \texttt{edge2element} is as follows:
\begin{lstlisting}[language=Matlab, caption={Vectorization of the \texttt{edge2element} matrix.},label=lst:vec_ed2el, escapechar=|]
%% to generate an element of edge
Inode=elements'; Inode=Inode(:); |\label{line:Inode}|% initial nodes
Enode=elements(:,[2,3,1])'; Enode=Enode(:);|\label{line:Enode}| % end nodes
Q=[Inode, Enode, nodes2element(sub2ind([nC nC],Inode,Enode)),nodes2element(sub2ind([nC nC],Enode,Inode))]; |\label{line:Q}|
p2 = nodes2edge(sub2ind([nC nC],Inode,Enode)); |\label{line:p2}|
[uniquep2,uniqueIdx] = unique(p2); |\label{line:ed2elStart}|% find the indices of the unique p2
r=setdiff(1:numel(p2),uniqueIdx)';  % find the indices of the repeated values p2 or find the repeated edges row indices of Q
Q(r,:)=[]; % removing the repeated edges in Q
edge2element(uniquep2,:)=Q; |\label{line:ed2elEnd}|
%% To produce the interior edges 
intE=find(edge2element(:,4));
%% To produce the exterior edges 
extE=find(edge2element(:,4)==0);
%% To produce local edge number
localE=repmat([3;1;2],nE,1); |\label{line:localE}|
ledTN(p2,1)=localE;  |\label{line:ledTN}|  % local edge numbering to T- elements
localE(r)=[]; % removing the repeated local edges
led(uniquep2,1)=localE; |\label{line:led}|% local edge numbering to T+ elements
\end{lstlisting}
In the above code
\begin{itemize}
\item Line \ref{line:Inode}-\ref{line:Enode}: The column vectors, \texttt{Inode} and \texttt{Enode}, specify the initial and end nodes of every edge $E \in \mathcal{E}^h$ through their respective indices.
\item Line \ref{line:Q}: The vector \texttt{Q} concatenates the \texttt{Inode}, \texttt{Enode}, $T_{+}$, and $T_{-}$, potentially resulting in the repetition of edges.
\item Line \ref{line:p2}:  \texttt{p2} is a column vector containing the global edge number according to $\text{conv} \{~\texttt{Inode},~\texttt{Enode}~\}$.
\item Line \ref{line:ed2elStart}-\ref{line:ed2elEnd}: We assemble the \texttt{edge2element} with \texttt{Q} by removing the repetition of edges.

\item Line \ref{line:localE}: The local edge numbering of element $T$ is $\{E_3,~ E_1,~ E_2\}$, that is $[3;1;2]$ in MATLAB.
\item Line \ref{line:ledTN}: With the help of repeated edges, we calculate the local edge numbering for $T_{-}$ elements in \texttt{ledTN}.
\item Line \ref{line:led}: With the help of unique edges, we compute the local edge numbering for $T_{+}$ elements in \texttt{led}.
\end{itemize}

\subsection{Vectorization of $\mathbf{redrefine.m}$}
The function  $\mathbf{redrefine\_vec.m}$ performs uniform refinement of the computational domain. 

\begin{lstlisting}[language=Matlab, caption={Vectorization of redrefine function.},label=lst:vec_redrefine, escapechar=|]
function [coordinates,elements,Dirichlet,Neumann]=redrefine_vec(coordinates,elements,n2ed,ed2el,Dirichlet,Neumann)
nE=size(elements,1);  % number of elements
nrEdges=size(ed2el,1);  % number of edges
%% Coordinates  |\label{line:CoordinateStart}|
Inode=ed2el(:,1); Enode=ed2el(:,2);  % initial and end nodes of edge
c1=coordinates(Inode,:); c2=coordinates(Enode,:); % their coordinates 
NCoord=(c1+c2)/2;  % coordinates of new nodes or mid-point of edges
Marker=(size(coordinates,1)+1:size(coordinates,1)+nrEdges);
coordinates(Marker,:)=[NCoord(:,1) NCoord(:,2)]; |\label{line:CoordinateEnd}|
%% Elements |\label{line:ElementsStart}|
CT=elements; % current triangle
CE= diag(n2ed(elements(:,[2 3 1]),elements(:,[3 1 2])));
M1=Marker(CE(1:nE)); M2=Marker(CE(nE+1:2*nE)); M3=Marker(CE(2*nE+1:3*nE));
elements((1:nE),:)=[M1' M2' M3'];
elements((nE+1:3:4*nE-2),:)=[CT(:,1) M3' M2'];
elements((nE+2:3:4*nE-1),:)=[CT(:,2) M1' M3'];
elements((nE+3:3:4*nE),:)=[CT(:,3) M2' M1'];  |\label{line:ElementsEnd}|
%% Dirichlet and Neumann boundary |\label{line:DirichletNeumannStartRefine}|
if (~isempty(Dirichlet)), Dirichlet=refineEdges(Dirichlet); end
if (~isempty(Neumann)), Neumann=refineEdges(Neumann); end
    function ed2n=refineEdges(ed2n)
        nEdges=size(ed2n,1);  % number of edges
        Base=diag(n2ed(ed2n(:,1),ed2n(:,2)));
        P=[ed2n(:,1) Marker(Base(:))' ed2n(:,2)];
        ed2n(1:nEdges,:)=[P(:,1) P(:,2)];
        ed2n(nEdges+1:2*nEdges,:)=[P(:,2) P(:,3)];
    end |\label{line:DirichletNeumannEndRefine}|
end
\end{lstlisting}

\begin{itemize}
\item Line \ref{line:CoordinateStart}-\ref{line:CoordinateEnd}: 
 The index of all new nodes is stored in a vector called \texttt{Marker}. The new midpoint coordinates are added to \texttt{coordinates} data by using \texttt{Marker} and concatenation.


\item Line \ref{line:ElementsStart}-\ref{line:ElementsEnd}: For all triangles in the original mesh, three new triangles are created by the mid-point of the edges. Mark the new nodal numbers in \texttt{M1}, \texttt{M2}, and \texttt{M3} vectors. The connectivity of all these new triangles is stored in the \texttt{elements} data.
\item  Line \ref{line:DirichletNeumannStartRefine}-\ref{line:DirichletNeumannEndRefine}: If there are $\mathcal{E}^h_{D}$ in the original mesh, collect the Dirichlet edge data in \texttt{Base}. A new node is created at the mid-point of the Dirichlet edge. Then for each Dirichlet edge, we find the corresponding new mid-point node using \texttt{Marker}. Update the \texttt{Dirichlet} data with these new nodes. For the Neumann boundary, the code does the same thing as for the Dirichlet boundary, but with the \texttt{Neumann} data.

\end{itemize}

\subsection{Run-time  of $\mathbf{edge\_vec.m}$ and $\mathbf{redrefine\_vec.m}$}
Table \ref{tab:redrefine_edge_table} presents the run-time of vectorized versions of the functions $\mathbf{edge.m}$ and $\mathbf{redrefine.m}$ in MATLAB, as we increase the level of refinement. Level $0$ corresponds to the initial triangulation of the domain (\figurename{~\ref{InitialTriangulation}}), which consists of four triangular elements. In contrast, level $1$ is the first uniform refinement of level $0$, resulting in 16 triangular elements. 
For any $n \in \mathbb{N}$, level $n$ is the uniform refinement of level $n-1$, and it contains $4^{n+1}$ triangular elements. 

\begin{table}[h]
\footnotesize
\captionsetup{font=small}
\vspace{-0.1cm}
\begin{tabular*}{\textwidth}{@{\extracolsep\fill}crll}
\vspace{-0.1cm}
Level & \thead{\texttt{nE}}   & \multicolumn{2}{@{}c@{}}{\thead{\textnormal{run-time (in sec)} \\  }} \\\cmidrule{3-4}%
 & & $\mathbf{edge\_vec.m}$ & $\mathbf{redrefine\_vec.m}$  \\
\midrule
5& $4096$  &    0.006 & 0.0001 \\  
6& $16384$ &   0.023 & 0.0008 \\   
7& $65536$  &   0.082   & 0.003\\  
8& $262144$  &  0.41 & 0.15 \\ 
9& $1048576$ &  2.04 &  0.68\\ 
10 &  $4194304$  & 8.66 & 2.90 \\
\end{tabular*}
\caption{Run-times: Vectorized code of $\mathbf{edge.m}$ and $\mathbf{redrefine.m}$.}
   \label{tab:redrefine_edge_table}
\end{table}

\begin{remark}
\vspace{-0.4cm}
Run-times for non-vectorized codes are much higher compared to vectorized codes; for instance, the run-times of $\mathbf{edge.m}$ and $\mathbf{redrefine.m}$ are 17114.35  and 2861.10 seconds, respectively, at level 9.
\end{remark}


\section{Implementation  of  second order elliptic model problem}\label{sec:ellipticproblem}
\subsection{Problem formulation}
Consider the general second-order elliptic problem
\begin{align}
-\nabla\cdot(A\nabla u+up)+\delta u&=f~\text{in}~~\Omega, \label{eq1}\\
~~u&=u_D~~\text{on}~~\Gamma_D,\label{eq2}\\
~~(A\nabla u+up)\cdot\nu&=g~~\text{on}~~\Gamma_N,\label{eq3}
\end{align}
where $A$ is a symmetric and positive definite matrix whose entries are constants, $p$ is a constant vector and $\delta$ is a constant. Also, $f,u_D$, and $g$ are assumed to be sufficiently smooth.

We consider the broken Sobolev space (cf. \cite{raviart}), i.e., a space of functions that belong to Sobolev space $H^{1}(T)$ on each element $T$,
$$
X=\{v\in L^{2}(\Omega): {v|}_{T}\in H^{1}(T),\forall T\in\mathcal{T}_{h}\}=\displaystyle\prod_{T\in\mathcal{T}_{h}}H^{1}(T).
$$
The broken Sobolev space has the norm 
$||v||_X=\bigg(\displaystyle\sum_{T\in\mathcal{T}_h}||\nabla v||^2_{L^2(T)}\bigg)^{\frac{1}{2}}.$
We shall denote by $H^{1/2}(\Gamma)$ the space of traces ${v|}_{\Gamma}$
over $\Gamma$ of the functions $v\in H^{1}(\Omega)$ and  $H^{-1/2}(\Gamma)$ denotes the dual space of $H^{1/2}(\Gamma)$.
 The Lagrange multiplier space is defined as
\begin{align*}
M=\bigg\{\chi\in\prod_{T\in\mathcal{T}_{h}}H^{-1/2}(\partial T/\Gamma_N):&~\text{there exists}~\textbf{q}\in H(\text{div},\Omega)
~\text{such that},~\textbf{q}\cdot\nu_{T}=\chi ~\text{on} ~\partial T/\Gamma_N, \forall T\in\mathcal{T}_{h}\bigg\}.
\end{align*}
The Lagrange multiplier space is provided with the norm 
$$\|\chi\|_M = \displaystyle\sup_{v\in X}\frac{\displaystyle\sum_{E\in\mathcal{E}^{h}/\mathcal{E}^h_{N}}\langle\chi ,v\rangle_E}{{\|v\|}_{X}},$$
where $\langle\cdot,\cdot\rangle_{E}$ denotes the duality paring between $H^{-1/2}(E)$ and
$H^{1/2}(E)$.

 Primal hybrid formulation of  \eqref{eq1}-\eqref{eq3} is to seek a pair of solutions
$(u,\kappa)\in X\times M$ such that
\begin{align}
a(u,v)+b(v, \kappa)&=(f,v)+\langle g,v\rangle_{\Gamma_N}\hspace{0.2cm}\forall v\in X,\label{eq1hyb1}\\
b(u,\chi)\qquad \qquad &=\langle u_D,\chi\rangle,\hspace{0.7cm}\forall \chi\in M,\label{eq1hyb2}
\end{align}
where 
\begin{align*}
a(u,v)&=\sum_{T\in\mathcal{T}_{h}}\int_{T}\left( (A\nabla u+up)\cdot \nabla v+\delta uv \right) \dx, \\
~b(v,\chi)&=-\sum_{E\in\mathcal{E}^{h}/\mathcal{E}^h_{N}}\langle\chi ,v\rangle_E,\\
\langle u_D,\chi\rangle&=-\sum_{E\in\mathcal{E}^{h}/\mathcal{E}^h_{N}}\langle\chi ,u_D\rangle_E.
\end{align*}
Here, the Lagrange multiplier
\begin{equation}\label{eq:exactLagrange}
\kappa=(A\nabla u+up)\cdot\nu_T\hspace{0.2cm}on~ \partial T/\Gamma_N,~ \forall T\in\mathcal{T}_{h}
\end{equation}
is associated with the constraint $u\in \left\{v\in H^{1}(\Omega):v=u_D~\text{on}~\Gamma_D\right\}$.
\subsection{Discrete problem and algebraic formulation}

Let $P_k$ be the space of polynomials of degree $\leq k$, for integer $k\geq 0$. We define the finite-dimensional spaces $X_h\subset X, M_h\subset M$ as 
\begin{align*}
&X_h=\prod_{T\in\mathcal{T}_h}P_1(T), \\
&M_{h}=\bigg\{\chi_{h}\in\prod_{E\in\mathcal{E}^h/\mathcal{E}^h_{N}}P_0(E):{\chi_{h}|}_{\partial T_{+}}
+{\chi_{h}|}_{\partial T_{-}}=0~\text{on}~T_{+}\cap T_{-},~
\text{for any adjacent pair}~T_{\pm}\in\mathcal{T}_{h}\bigg\}.
\end{align*}


We seek a pair $(u_h,\kappa_h)\in X_h\times M_h$ such that for all pairs $(v_h,\chi_h)\in X_h\times M_h$
\begin{align}
\sum_{T\in\mathcal{T}_h}\int_{T}((A\nabla u_h+ u_hp)\cdot\nabla  v_h+\delta u_hv_h)\dx&-\sum_{E\in\mathcal{E}^h/\mathcal{E}^h_{N}}\int_{E}\kappa_h[\![v_h]\!]_E\dgamma\nonumber\\
&=\sum_{T\in\mathcal{T}_h}\int_Tfv_h\dx+\int_{\Gamma_N}gv_h\dgamma,\label{ph21}\\
-\sum_{E\in\mathcal{E}^h/\mathcal{E}^h_{N}}\int_{E}\chi_h[\![u_h]\!]_E \dgamma&=-\sum_{E\in\mathcal{E}^h/\mathcal{E}^h_{N}}\int_{E}\chi_hu_{Dh}\dgamma.\label{ph22}
\end{align}
Note that the problem \eqref{ph21}-\eqref{ph22} has a unique pair of solution $(u_h,\kappa_h)\in X_h\times M_h$ (cf. \cite{raviart}).
With  $N=\dim(X_h)$, let $X_h=\spn\{\phi_1,\dots,\phi_N\}$ and $U=[x_1,\dots,x_N]^\prime$ are the components of $$u_h=\displaystyle\sum^N_{i=1}x_i\phi_i\in X_h.$$

Let $E_1, E_2, E_3$ be the edges of the triangle $T$,
and let $\nu_{E_k}$ denote the unit normal vector of $E_k$ chosen with a global fixed orientation while $\nu_k$ denotes the outer unit normal of $T$ along $E_k$.
We define the basis functions $\psi_E$ of $M_h$, where $E\in\mathcal{E}^h/\mathcal{E}^h_N$ as below:
\begin{equation*}
\psi_{E_k}(x) = \sigma_k ~~~ \mathrm{for}~ k=1,2,3~ \mathrm{and}~x \in T,
\end{equation*}
where $\sigma_k = \nu_{E_k}.\nu_{k}$ is $+1$ if $\nu_{E_k}$ points outward and otherwise~ $-1$.
Globally, let for any edge $E=T_+\cap T_-\in\mathcal{E}^h_\Omega$
\begin{align*}
\psi_E(x)=\begin{cases}
& 1~~~~~\text{for}~~x\in T_{+},\\
&- 1~~~\text{for}~~x\in T_{-},\\
&0~~~~~~\text{elsewhere};
\end{cases}
\end{align*}
and for $E\in \mathcal{E}^h_D$
\begin{align*}
\psi_E(x)=\begin{cases}
&1~~~~~\text{for}~~x\in E,\\
&0~~~~~\text{elsewhere}.
\end{cases}
\end{align*}

 With $L=\dim(M_h)$, let $M_h= \spn\{\psi_l\}_1^L$ with $\psi_l=\psi_{E_l}$, where $l=1,\dots, L$ is an enumeration of edges $$\mathcal{E}^h_\Omega\cup\mathcal{E}^h_{D}=\{E_{1},\dots,E_{L}\}.$$
 Let $\Lambda=[x_{N+1},\dots,x_{N+L}]^\prime$ are the components of $\kappa_h$ such that $$\kappa_h=\displaystyle\sum^{L}_{l=1}x_{N+l}\psi_l.$$ The problem \eqref{ph21}-\eqref{ph22} can be written in a linear system of equations for unknowns $U$ and $\Lambda$ as: for $j=1,\dots,N$ and $l=1,\dots,L$.
\begin{align}
\sum^N_{i=1}x_i\sum_{T\in\mathcal{T}_h}\int_T((A\nabla\phi_j+\phi_jp)\cdot\nabla\phi_i+\delta\phi_j\phi_i)\dx&-\sum_{l=1}^{L}x_{N+l}\int_{E_{l}}\psi_l[\![\phi_j]\!]_{E_l}\dgamma\nonumber\\
&=\sum_{T\in\mathcal{T}_h}\int_{T}f\phi_j\dx+\int_{\Gamma_N}g\phi_j\dgamma,\label{ph23}\\
-\sum_{i=1}^{N}x_{i}\int_{E_l}\psi_l[\![\phi_i]\!]_{E_l}\dgamma&= - \int_{E_l}\psi_lu_{Dh}\dgamma.\label{ph24}
\end{align}

\subsection{Local matrices $\, \mathbb{B}_T$, $\mathbb{D}_T$ $\mathbb{C}_T$, 
$\mathbb{M}_T$}
With local numbers $i\in\{1,2,3\}$,  let $\{a_i\}^{3}_{i=1}$ be the set of vertices  of an element $T\in\mathcal{T}_h$ and  $\{\lambda_{i}(x)\}_{i=1}^{3}$ be the set of barycentric coordinates of a point $x\in T$ (cf. \cite{ciar}). Let   $\{\phi^{T}_{i}(x)\}^{3}_{i=1}$
 be the set of basis functions of $P_1(T)$.
 It should be noted that the basis functions of $X_h$ are not globally continuous over $\Omega$ but locally continuous, which means that for any triangle $T\in \mathcal{T}_h$ with global node numbers $i, j, k$, $\phi_i=\lambda^T_1; \phi_j=\lambda^T_2; \phi_k=\lambda^T_3$. \\
 
 \begin{defn}\label{def:localMatrices}
 The local matrices 
 are defined  for $i,j,l=1,2,3$ as follows:
\begin{align*}
&(\mathbb{B}_T)_{ji}= \int_T A \nabla \phi_j \cdot \nabla \phi_i ~\dx, \qquad 
(\mathbb{D}_T)_{ji}= \int_T  (\phi_j p) \cdot \nabla \phi_i \dx, \\
&(\mathbb{C}_T)_{lj}= \int_{E_{l}}\psi_{l}[\![\phi_j]\!]_{E_l}\dgamma, \qquad 
(\mathbb{M}_T)_{ji}= \int_T  \delta \phi_j  \phi_i \dx.
\end{align*}
The elaborate computation provides forms of all local matrices:
\begin{align}\label{local_B}
\mathbb{B}_T & =|T|\left(
\setlength{\tabcolsep}{20pt}
\renewcommand\arraystretch{2}
    \begin{array}{ccc}
      (A\nabla\lambda^T_1)\cdot\nabla\lambda^T_1 & (A\nabla\lambda^T_1)\cdot\nabla\lambda^T_2& (A\nabla\lambda^T_1)\cdot\nabla\lambda^T_3 \\
    (A \nabla\lambda^T_2)\cdot\nabla\lambda^T_1 & (A\nabla\lambda^T_2)\cdot\nabla\lambda^T_2 & (A\nabla\lambda^T_2)\cdot\nabla\lambda^T_3\\
     (A\nabla\lambda^T_3)\cdot\nabla\lambda^T_1 & (A\nabla\lambda^T_3)\cdot\nabla\lambda^T_2 & (A\nabla\lambda^T_3)\cdot\nabla\lambda^T_3
    \end{array}
    \right), \\
\label{local_D}
\mathbb{D}_T & =\frac{|T|}{3}\left(
\setlength{\tabcolsep}{20pt}
\renewcommand\arraystretch{2}
    \begin{array}{ccc}
      p\cdot\nabla\lambda^T_1 & p\cdot \nabla\lambda^T_1& p\cdot\nabla\lambda^T_1\\
     p\cdot\nabla\lambda^T_2 & p\cdot\nabla\lambda^T_2& p\cdot\nabla\lambda^T_2\\
     p\cdot\nabla\lambda^T_3 & p\cdot\nabla\lambda^T_3& p\cdot\nabla\lambda^T_3
    \end{array}
    \right), \\
\label{local_M}
\mathbb{M}_T & = \delta \left(
\setlength{\tabcolsep}{20pt}
\renewcommand\arraystretch{2}
    \begin{array}{ccc}
      \displaystyle\int_{T}\lambda^T_1\lambda^T_1\dx & \displaystyle\int_{T}\lambda^T_1\lambda^T_2\dx& \displaystyle\int_{T}\lambda^T_1\lambda^T_3\dx \\
     \displaystyle \int_{T}\lambda^T_2\lambda^T_1\dx & \displaystyle\int_{T}\lambda^T_2\lambda^T_2\dx & \displaystyle\int_{T}\lambda^T_2\lambda^T_3\dx\\
     \displaystyle \int_{T}\lambda^T_3\lambda^T_1\dx & \displaystyle\int_{T}\lambda^T_3\lambda^T_2\dx & \displaystyle\int_{T}\lambda^T_3\lambda^T_3\dx
    \end{array}
    \right)=\delta  |T|\left(
    \setlength{\tabcolsep}{20pt}
\renewcommand\arraystretch{2}
    \begin{array}{ccc}
      \frac{1}{6} & \frac{1}{12}& \frac{1}{12} \\
     \frac{1}{12}  & \frac{1}{6} & \frac{1}{12} \\
     \frac{1}{12}  & \frac{1}{12}  & \frac{1}{6}
    \end{array}
    \right), \\
\label{local_C}
\mathbb{C}_T & =\left(
\setlength{\tabcolsep}{20pt}
\renewcommand\arraystretch{2}
    \begin{array}{ccc}
      \displaystyle\int_{E_1}\psi_1[\![\lambda^T_1]\!]_{E_1}\dgamma & \displaystyle\int_{E_1}\psi_1[\![\lambda^T_2]\!]_{E_1}\dgamma& \displaystyle\int_{E_1}\psi_1[\![\lambda^T_3]\!]_{E_1}\dgamma\\
     \displaystyle\int_{E_2}\psi_2[\![\lambda^T_1]\!]_{E_2}\dgamma & \displaystyle\int_{E_2}\psi_2[\![\lambda^T_2]\!]_{E_2}\dgamma & \displaystyle\int_{E_2}\psi_2[\![\lambda^T_3]\!]_{E_2}\dgamma\\
     \displaystyle\int_{E_3}\psi_3[\![\lambda^T_1]\!]_{E_3}\dgamma & \displaystyle\int_{E_3}\psi_3[\![\lambda^T_2]\!]_{E_3}\dgamma & \displaystyle\int_{E_3}\psi_3[\![\lambda^T_3]\!]_{E_3}\dgamma
    \end{array}
    \right)=\left(
    \setlength{\tabcolsep}{20pt}
\renewcommand\arraystretch{2}
    \begin{array}{c}
       |E_1|\sigma_1 R_1  \\
      |E_2| \sigma_2 R_2  \\
     |E_3| \sigma_3 R_3
    \end{array}
  \right),
\end{align}
where
$
R_1=\left[
    \begin{array}{c}
      0,  ~~
       \frac{1}{2},
       ~~\frac{1}{2}
    \end{array}
  \right],\quad  R_2=\left[
    \begin{array}{c}
      \frac{1}{2},  ~~
       0,
       ~~\frac{1}{2}
    \end{array}
  \right], \quad R_3=\left[
    \begin{array}{c}
      \frac{1}{2},  ~~
       \frac{1}{2},
       ~~0
    \end{array}
  \right].
$
\end{defn}
\subsection{Right-hand side and boundary conditions}\label{sec:righthand} 
The computation of the right-hand side in \eqref{ph23}-\eqref{ph24} includes numerical integration over elements and edges of the functions $f,~g,$ and $u_D$. The right-hand side $F=\{F_j\}_{j = 1}^{N}$ in \eqref{ph23} is defined as
$
F_j
= b_j+LN_j ~, 
$
where
$
b_j= \displaystyle \sum_{T\in\mathcal{T}_h} \int_Tf\phi_j \dx~~\mathrm{and}~~LN_j= \displaystyle \sum_{E\in\mathcal{E}_N^h} \int_{E}g\phi_j \dgamma.
$
Both parts are approximated by
 \begin{align}\label{eq:right_b}
& b_j 
\approx \displaystyle \sum_{T\in\mathcal{T}_h} \frac{|T|}{3}~\sum_{i=1}^{3} f(m_{E_{i}}) \, \phi_j(m_{E_{i}}), \\
\label{eq:right g}
& LN_j \approx \displaystyle \sum_{E\in\mathcal{E}_N^h} |E| g(m_{E}) \, \phi_j(m_{E}),~~~~~~\text{where}~ m_{E}\in \mathcal{N}_m(\Gamma_N).
 \end{align}
The right hand side $b_D=\{b_{D_l}\}_{l = 1}^{L}$ in \eqref{ph24} is defined and approximated as
\begin{align}\label{eq:DirichletEqn}
b_{D_l}&=-\displaystyle\int_{E_l} \psi_lu_{Dh}\dgamma \nonumber\\ &\approx -  |E_l| \, \sigma_l u_{Dh} \,(m_{E_l}),~~~~\text{where}~m_{E_l}\in \mathcal{N}_m(\Gamma_D) \nonumber\\ & \approx -|E_l| \, u_{Dh}(m_{E_l}), ~~~~~\text{since}~\sigma_l=1~\text{for}~l=1,\dots,L.
\end{align}

\subsection{Assembly of the global matrices}
In vector-matrix form, \eqref{ph23}-\eqref{ph24} can be written as


\begin{equation}\label{globalSystem}
\mathrm{M}_{\mathrm{vec}}\mathrm{W}_{\mathrm{vec}}=\mathrm{F}_{\mathrm{vec}}~,
\end{equation}
where global matrices and vectors read
\begin{align*}
&\mathrm{M}_{\mathrm{vec}} =\left(
    \begin{array}{cc}
      \mathbb{B}+\mathbb{D}+\mathbb{M} & -\mathbb{C}^{\prime}  \\
      -\mathbb{C} & 0  \\
    \end{array}
  \right)_{(N+L)\times (N+L)},   
  \qquad 
  \mathrm{W}_{\mathrm{vec}}= \left(\begin{array}{c}
      U   \\
      \Lambda  \\
    \end{array}
  \right)_{(N+L)\times 1}, 
  \quad 
  \mathrm{F}_{\mathrm{vec}}=\left(
    \begin{array}{c}
      F   \\
      b_D  \\
    \end{array}
  \right)_{(N+L)\times 1}, \\
&\mathbb{B}=\big(\mathbb{B}_{ij}\big)_{N\times N}, \quad \mathbb{B}_{ij}=\displaystyle\sum_{T\in\mathcal{T}_h}\int_{T}(A\nabla \phi_i)\cdot\nabla\phi_j~\dx, \\
&\mathbb{D}=\big(\mathbb{D}_{ij}\big)_{N\times N}, \quad 
\mathbb{D}_{ij}
=\displaystyle\sum_{T\in\mathcal{T}_h}\int_{T}  (\phi_i p)\cdot\nabla\phi_j~\dx, \\
&\mathbb{M}=\big(\mathbb{M}_{ij}\big)_{N\times N}, \quad
\mathbb{M}_{ij}
=\displaystyle\sum_{T\in\mathcal{T}_h}\int_{T}  \phi_i \cdot \phi_j~\dx, \\
&\mathbb{C}=(\mathbb{C}_{N+l,i})_{L\times N}, \quad 
\mathbb{C}_{N+l,i}=\displaystyle \int_{E_{l}}\psi_l[\![\phi_i]\!]_{E_l}\dgamma,\\
&F=(F_j)_{ N\times 1}, \quad b_D=(b_{D_l})_{L\times 1}.
\end{align*}

\subsection{Implementation}

\subsubsection*{MATLAB code of vectorized $\mathbf{StiffMassConv\_PH.m}$}

Here we  explain the  vectorized code \textbf{StiffMassConv\_PH.m}, where we compute the matrices $\mathbb{B},~ \mathbb{D},~ \mathbb{M} \in \mathbb{R}^{N \times N}$, and $b\in \mathbb{R}^{N \times 1}$. The vectorization of stiffness matrix $\mathbb{B}$ and load vector $b$ are similar to the MATLAB code Listing 3 in  \cite{funken2011efficient} with an adequate amount of modifications.
\begin{lstlisting}[language=Matlab, caption= \small{Vectorization of the assembly of stiffness matrix $\mathbb{B}$, convection matrix $\mathbb{D}$, mass matrix $\mathbb{M}$ and load vector $b$.}, label=lst:vec_StiffMassConv,escapechar=|]
function [B,D,M,b]=StiffMassConv_PH(coordinates,elements,params)|\label{line:StiffMassConvPH}|
A=params.A;p=params.p;delta=params.delta;
nE=size(elements,1); |\label{line:numberofElements}| % number of elements
%% Vectorization of B matrix
P1=coordinates(elements(:,1),:); % collection of the first nodes of all triangles
P2=coordinates(elements(:,2),:); % collection of the second nodes of all triangles
P3=coordinates(elements(:,3),:); % collection of the third nodes of all triangles
d21=P2-P1;  d31=P3-P1;  d32=P3-P2;  % represents the edges of the elements in a mesh
mp23=1/2.*(P2+P3);  mp13=1/2.*(P1+P3);   mp12=1/2.*(P1+P2); % midpoints of the edges
area4=2*(d21(:,1).*d31(:,2)-d21(:,2).*d31(:,1));  % 4*Area
GI=reshape([1:3*nE; 1:3*nE; 1:3*nE],9,nE); |\label{line:GI}| % 3*nE =Number of nodes
I3=reshape(GI([1 4 7 2 5 8 3 6 9],:),9*nE,1); 
J3=reshape(GI,9*nE,1); |\label{line:J3}|
a3=(sum(([d21(:,2),d21(:,1)]*A).*[d31(:,2),d31(:,1)],2)./area4)'; |\label{line:a3}|
b3=(sum(([d31(:,2),d31(:,1)]*A).*[d31(:,2),d31(:,1)],2)./area4)';
c3=(sum(([d21(:,2),d21(:,1)]*A).*[d21(:,2),d21(:,1)],2)./area4)';
BT=[-2*a3+b3+c3; a3-b3; a3-c3; a3-b3;b3;-a3;a3-c3;-a3;c3];
B=sparse(I3,J3,BT,3*nE,3*nE); |\label{line:B}|
%% Vectorization of mass matrix (M) |\label{line:MassStart}|
MT=reshape(delta.*(area4'./48).*[2;1;1;1;2;1;1;1;2],9*nE,1);
M=sparse(I3,J3,MT,3*nE,3*nE); |\label{line:M}|
%% Vectorization of D matrix  |\label{line:ConvectionStart}|
h1=(-p(1)*d32(:,2)+p(2)*d32(:,1))'./6;
h2=(p(1)*d31(:,2)-p(2)*d31(:,1))'./6;
h3=(-p(1)*d21(:,2)+p(2)*d21(:,1))'./6;
DT=[ h1; h1; h1; h2; h2; h2; h3; h3; h3]; 
D=sparse(I3,J3,DT,3*nE,3*nE); |\label{line:D}|
%% Vectorization of load vector b |\label{line:loadStart}|
f1=f(mp13,params)./2 +f(mp12,params)./2;
f2=f(mp23,params)./2+f(mp12,params)./2;
f3=f(mp13,params)./2+f(mp23,params)./2;
fT=(area4.*[f1 f2 f3])./12;
b = accumarray (reshape(1:3*nE,3*nE,1), reshape(fT',3*nE,1) ,[3*nE 1]); |\label{line:loadb}|
end
\end{lstlisting}

\begin{itemize}
\item Line \ref{line:GI}-\ref{line:J3}: We create a $9 \times \texttt{nE}$ matrix called \texttt{GI} by reshaping the matrix $[1:3 \times \texttt{nE}; 1:3 \times \texttt{nE}; 1:3 \times \texttt{nE}]$ such that each column of \texttt{GI} (cf. \figurename{~\ref{fig:GIdiagram}}) corresponds to the column index of local stiffness matrices for elements in the global assembly. \texttt{I3} and \texttt{J3} (cf. \figurename{~\ref{fig:stiffPH_diagram}}) contain index pairs for the global stiffness matrix, constructed using \texttt{GI}. Note that, \texttt{I3} is based on the modified \texttt{GI} with specific row indices $[1~ 4~ 7~ 2 ~5 ~8~ 3~ 6~ 9]$, while \texttt{J3} is based on the original \texttt{GI}.
\newline Define, the row and column index pair of entries in the stiffness matrix corresponding to the $m^{\text{th}}$ element as 
$$ \rho(m) = \begin{bmatrix}
            3(m-1)+1 \\
            3(m-1)+2\\
            3(m-1)+3\\
            3(m-1)+1\\
            3(m-1)+2\\
            3(m-1)+3\\
            3(m-1)+1\\
            3(m-1)+2\\
            3(m-1)+3
            \end{bmatrix}
 ~~\text{and}~~
 \theta(m) = \begin{bmatrix}
            3(m-1)+1 \\
            3(m-1)+1\\
            3(m-1)+1\\
            3(m-1)+2\\
            3(m-1)+2\\
            3(m-1)+2\\
            3(m-1)+3\\
            3(m-1)+3\\
            3(m-1)+3
            \end{bmatrix},~~~ \text{where} ~m=1,2,\dots,\texttt{nE}. $$
           
\begin{figure}[H]
\centering
\includegraphics[width=7.5cm]{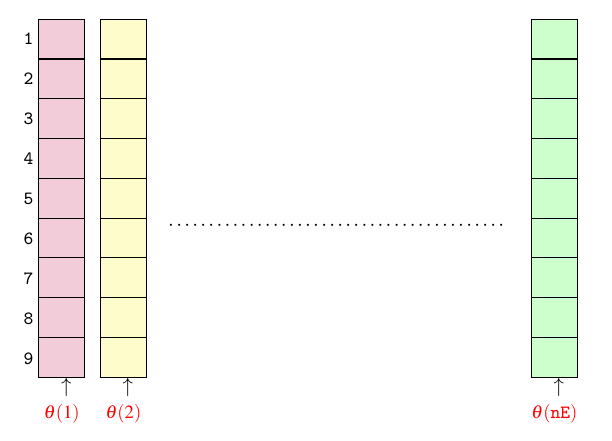}
\centering
\caption{Geometrical representation of the columns of the \texttt{GI} matrix in MATLAB.}
\label{fig:GIdiagram}
 \end{figure}            
\begin{figure}[H]
\centering
\includegraphics[width=13.5cm]{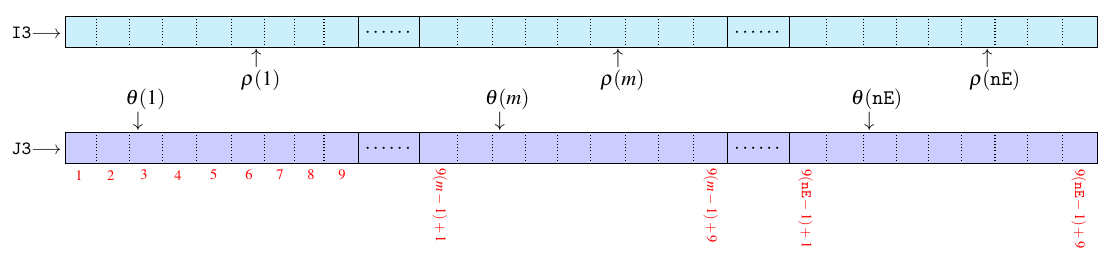}
\centering
\caption{Geometrical representation of \texttt{I3} and \texttt{J3} column vectors in MATLAB.}
\label{fig:stiffPH_diagram}
\end{figure} 

 \begin{remark}
The index pair vectors for conforming (cf. \cite{funken2011efficient}) and non-conforming FEM differ in vectorization because of the distinct element-wise assembly of the local stiffness matrices. For primal hybrid FEM, we generate \texttt{I3} and \texttt{J3} index pair vectors that are compatible with non-conforming P1-FEM.
\end{remark}
 \item Line \ref{line:a3}-\ref{line:B}: 
 The vectors \texttt{a3}, \texttt{b3}, and \texttt{c3} for element $T$ are
\begin{align}\label{eq:a3}
            {\texttt{a3}|}_{T} & =\frac{1}{4|T|} \left(  \begin{bmatrix}
            y_2-y_1 & x_2-x_1
            \end{bmatrix} A \right)
\cdot\begin{bmatrix}
            y_3-y_1 & x_3-x_1
            \end{bmatrix}, \\
\label{eq:b3}
           {\texttt{b3}|}_{T} & =\frac{1}{4|T|} \left(  \begin{bmatrix}
            y_3-y_1 & x_3-x_1
            \end{bmatrix} A \right)
\cdot \begin{bmatrix}
            y_3-y_1 & x_3-x_1
            \end{bmatrix}, \\
             {\texttt{c3}|}_{T} & =\frac{1}{4|T|} \left( \begin{bmatrix}
            y_2-y_1 & x_2-x_1
            \end{bmatrix} A \right)
\cdot \begin{bmatrix}
            y_2-y_1 & x_2-x_1
            \end{bmatrix}.
 \end{align}

\item Line \ref{line:MassStart}-\ref{line:M}:
The element-wise assembly of $\mathbb{M}$ uses nine updates of the mass matrix per element $T$. That is, the vector $\texttt{MT}$ has length $9 \times $ \texttt{nE}.

 \item Lines \ref{line:ConvectionStart}-\ref{line:D}:  ${\texttt{h1}|}_{T}$, ${\texttt{h2}|}_{T}$, and ${\texttt{h3}|}_{T}$ are  entries of the matrix $\mathbb{D}_T$ (\ref{local_D}),
     \begin{align*}
     {\texttt{h1}|}_{T} = \frac{1}{6}~ p\cdot \begin{bmatrix}
     -(y_3-y_2) \\
     x_3-x_2
     \end{bmatrix}, \quad 
     ~{\texttt{h2}|}_{T} = \frac{1}{6}~ p\cdot\begin{bmatrix}
     y_3-y_1 \\
     -(x_3-x_1)
     \end{bmatrix}, \quad 
     ~{\texttt{h3}|}_{T} = \frac{1}{6}~ p\cdot \begin{bmatrix}
     -(y_2-y_1) \\
     x_2-x_1
     \end{bmatrix}.
     \end{align*}
   We assemble the convection matrix $\mathbb{D}$ by using \textbf{sparse} command in MATLAB with \texttt{I3}(row index) and \texttt{J3}(column index).

\item Line \ref{line:loadStart}-\ref{line:loadb}: We evaluate the volume forces $f(m_{E_{i}}),~ i=1,2,3$ for all triangle $T \in \mathcal{T}_h$, then  \eqref{eq:right_b} and assemble the load vector $b$ using $\mathbf{accumarray}$.
\end{itemize}

\subsubsection*{MATLAB code of Vectorized $\mathbf{Lambda\_PH.m}$}

\begin{lstlisting}[language=Matlab, caption={Vectorized and efficient MATLAB implementation of the assembly of matrix $\mathbb{C}$.}, ,label=lst:vec_lambdaPH, escapechar=|]
 function [C,e11,v,R11]= Lambda_PH(coordinates,ed2el,intE,Dbed,Redges,led,ledTN) |\label{line:LambdaPHfun}|
nrEdges=size(ed2el,1); % number of edges
num_intE=size(intE,1); % number of interior edges
nDb=size(Dbed,1); % number of Dirichlet boundary edges
edges=ed2el(1:nrEdges,1:2);   % all edges 
%% Vectorization of matrix C
v=vecnorm((coordinates(edges(:,1),:)-coordinates(edges(:,2),:))')';  |\label{line:v}| % length of edges
e11=ed2el(:,3); e12=nonzeros(ed2el(:,4)); |\label{line:e11e12}|
i=find(led==1); j=find(led==2); k=find(led==3); |\label{line:ijk}|
p=find(ledTN(intE)==1); q=find(ledTN(intE)==2); r=find(ledTN(intE)==3); |\label{line:pqr}|
R1=[0, 1/2, 1/2]; R2= [1/2, 0,  1/2]; R3= [1/2, 1/2, 0]; |\label{line:R1R2R3}|
R11(i,:)=repmat(R1,size(i));
R11(j,:)=repmat(R2,size(j));
R11(k,:)=repmat(R3,size(k));
R22(p,:)=repmat(R1,size(p));
R22(q,:)=repmat(R2,size(q));
R22(r,:)=repmat(R3,size(r)); |\label{line:R22end}|
I1=reshape([Dbed' intE';Dbed' intE';Dbed' intE'],3*(num_intE+nDb),1); |\label{line:I1}|
J1=reshape(((e11([Dbed' intE'])-ones(num_intE+nDb,1))*3 +[1 2 3])'...
        ,3*(num_intE+nDb),1);
S1=reshape(( v([Dbed' intE']).* R11([Dbed' intE'],:))',3*(num_intE+nDb),1);
I2=reshape([intE';intE';intE'],3*num_intE,1);
J2=reshape(((e12(:)-ones(num_intE,1)).*3 +[1 2 3])',3*num_intE,1);
S2=reshape(( (-1).*v(intE).*(R22))',3*num_intE,1); |\label{line:S2}|
I=[I1;I2]; J=[J1;J2]; S=[S1;S2]; |\label{line:append}|% appending vectors
C=sparse(I, J, S);
C=C(Redges, :);  |\label{line:C}| % eliminate extra rows
end



\end{lstlisting}
\begin{itemize}
\item Line \ref{line:LambdaPHfun}: The function $\mathbf{Lambda\_PH.m}$ takes the \texttt{coordinates} data, \texttt{ed2el}, \texttt{intE}, Dirichlet boundary edges(\texttt{Dbed}),  \texttt{Redges}( $\mathcal{E}^h/\mathcal{E}^h_N$), \texttt{led},  and \texttt{ledTN}. It computes the  matrix $\mathbb{C} \in \mathbb{R}^{L \times N }$ in \eqref{globalSystem} and other required data $\texttt{e11},~\texttt{v},~\texttt{R11}$ for $\mathbf{Main.m}$ described in the next subsection.
\item Line \ref{line:v}-\ref{line:e11e12}:
The length of all edges is computed by using the MATLAB function $\mathbf{vecnorm}$ and stored in a vector \texttt{v}. Vector $\texttt{e11}=\texttt{ed2el}(:,3)$, which contains the index of all $T_{+}$ elements, and  vector $\texttt{e12}=\mathbf{nonzeros}(\texttt{ed2el}(:,4))$, which contains the index of all $T_{-}$ elements.
\item Line \ref{line:ijk}-\ref{line:pqr}: For each triangular element $T_{+}\in \mathcal{T}_h$, the vectors \texttt{i}, \texttt{j}, \texttt{k} store the indices of the first edge ($E_1$), second edge ($E_2$), third edge ($E_3$), respectively.   Similarly, the vectors \texttt{p}, \texttt{q}, \texttt{r} store the indices of $E_1$, $E_2$, $E_3$ for $T_{-}$ elements.

\item  Line \ref{line:R1R2R3}-\ref{line:R22end}: We construct the matrices $\texttt{R11} \in \mathbb{R}^{L \times 3}$ and $ \texttt{R22}\in \mathbb{R}^{N_{\Omega} \times 3}$ based on the values of the \texttt{led} and \texttt{ledTN} vectors, as well as the $ R_1,~R_2, ~\text{and}~ R_3$ vectors (Definition \ref{def:localMatrices}). The rows of \texttt{R11} corresponding to $E_1$, $E_2$, and $E_3$ of each $T_{+} \in \mathcal{T}_h$ are populated with the vectors $R_1$, $R_2$, and $R_3$, respectively. This is accomplished by using the $\mathbf{repmat}$ function to replicate each vector along the rows of \texttt{R11} at the indices stored in \texttt{i}, \texttt{j}, and \texttt{k}, respectively. Similarly, for each $T_{-} \in \mathcal{T}_h$,  we build the rows of \texttt{R22} are $R_1$, $R_2$, and $R_3$ vectors with respect to the indices stored in \texttt{p}, \texttt{q}, and \texttt{r}.

\item Line \ref{line:I1}-\ref{line:S2}: We calculate the row and column indexing vectors \texttt{I1} and \texttt{J1} $\in \mathbb{R}^{3L \times 1}$ (as shown in \figurename{~\ref{fig:LambdaPH_diagram1}}) for the matrix $\mathbb{C}$ based on the $T_{+}$ elements. Similarly, we obtain the row and column indexing vectors \texttt{I2} and \texttt{J2} $\in \mathbb{R}^{3N_{\Omega} \times 1}$ (as shown in \figurename{~\ref{fig:LambdaPH_diagram2}}) for the matrix $\mathbb{C}$ based on the $T_{-}$ elements. The entries of matrix $\mathbb{C}$ corresponding to  $T_{+}$ and $T_{-}$ elements are represented by \texttt{S1} $\in \mathbb{R}^{3L \times 1}$ (as shown in \figurename{~\ref{fig:LambdaPH_diagram1}}) and \texttt{S2} $\in \mathbb{R}^{3N_{\Omega} \times 1}$ (as shown in \figurename{~\ref{fig:LambdaPH_diagram2}}), respectively, and are arranged in column vectors:
\begin{align*}
\texttt{S1}_{|_{T_{+}}} &=  \begin{pmatrix} \displaystyle \texttt{v}(i) ~ \texttt{R11}(i,:)  \\
      \displaystyle \texttt{v}(j)~   \texttt{R11}(j,:)  \\
     \displaystyle \texttt{v}(k) ~ \texttt{R11}(k,:)  \end{pmatrix},~ \text{where}~ i,~j,~\text{and}~k ~\text{are the global edge number of} ~T_{+}, \\
\text{and}~~
\texttt{S2}_{|_{T_{-}}}&= \begin{pmatrix} -\texttt{v}(p) ~  \texttt{R22}(p,:)  \\
      -\texttt{v}(q) ~   \texttt{R22}(q,:)  \\
     -\texttt{v}(r) ~  \texttt{R22}(r,:)  \end{pmatrix} , ~ \text{where}~ p,~q,~\text{and}~r ~\text{are the global edge number of} ~T_{-}.
\end{align*}

Let $\alpha_1, \alpha_2, \dots, \alpha_{N_{D}}$ be the global numbers of the $\mathcal{E}_{D}^{h}$ for the mesh $\mathcal{T}_h$. Let $\beta_1, \beta_2, \dots, \beta_{N_{\Omega}}$ be global numbers of $\mathcal{E}_{\Omega}^{h}$ for the mesh $\mathcal{T}_h$.
Define,
\begin{align*}
 \gamma_m &= (\texttt{e11}(m)-1)*3+[1~2~3],~~~ \text{where} ~m=1,2,\dots,\texttt{nrEdges},  \\
\text{and}~~~~
\zeta_l & = (\texttt{e12}(l)-1)*3+[1~2~3],~~~ \text{where} ~l=1,2,\dots,N_{\Omega}. 
\end{align*}

 \begin{figure}
\centering
 \includegraphics[ width=13cm]{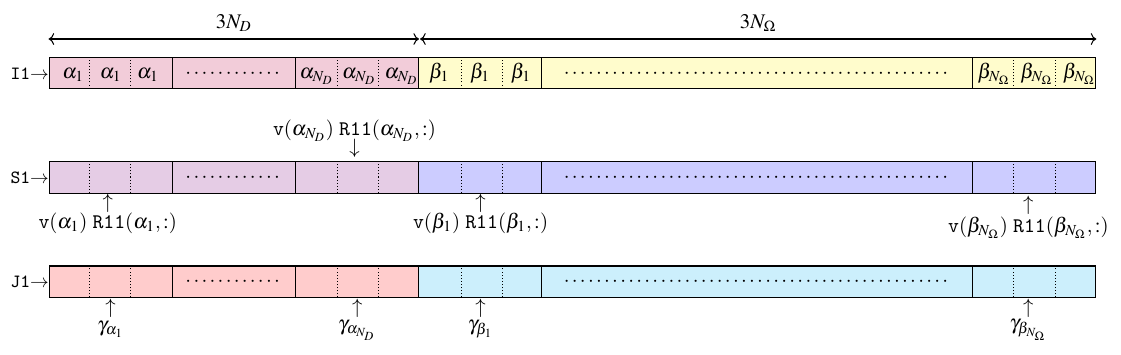}
\centering
\caption{Geometrical representation of \texttt{I1}, \texttt{J1}, and \texttt{S1} vectors in MATLAB.}
\label{fig:LambdaPH_diagram1}
\centering
\includegraphics[ width=10cm]{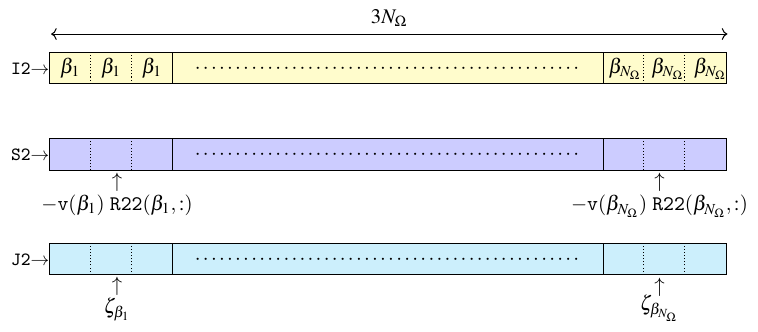}
\centering 
\caption{Geometrical representation of \texttt{I2}, \texttt{J2}, and \texttt{S2} vectors in MATLAB.}
\label{fig:LambdaPH_diagram2}
\end{figure}

\item Line \ref{line:append}-\ref{line:C}: We are appending the vectors \texttt{I=[I1;I2]}, \texttt{J=[J1;J2]}, \texttt{S=[S1;S2]}. We assemble the matrix $\mathbb{C}$ by using $\mathbf{sparse}$ with \texttt{I}, \texttt{J}, and \texttt{S}, but it has extra zero rows with respect to the $\mathcal{E}_{N}^{h}$. We eliminate the extra rows from matrix $\mathbb{C}$ in line \ref{line:C}.
\end{itemize}

\subsubsection*{MATLAB code of $\mathbf{Main.m}$}
\begin{lstlisting}[caption=$\mathbf{Main.m.}$, label=lst:Main,escapechar=|]
load('coordinates.dat');   load('elements.dat'); 
load('Dirichlet.dat');     load('Neumann.dat');
[n2el,n2ed,ed2el,intE,extE,led,ledTN]=edge_vec(elements,coordinates);
params.A=[1 0; 0 1]; params.p=[1,1]; params.delta=1; % problem parameters
nt=4;
 for p1=1:nt
    fprintf('Level = %d\n',p1); 
    tic  
    [coordinates,elements,Dirichlet,Neumann]=redrefine_vec(coordinates,elements,n2ed,ed2el,Dirichlet,Neumann);
    fprintf('run-time(in sec) for the function redrefine_vec.m  is = %d\n',toc);
    tic
    [n2el,n2ed,ed2el,intE,extE,led,ledTN]=edge_vec(elements,coordinates);
    fprintf('run-time(in sec) for the function edge_vec.m  is = %d\n',toc);
    h(p1)=sqrt(det([1 1 1;coordinates(elements(1,:),:)'])*2); % space parameter
    grad4e = getGrad4e(coordinates,elements);
    nrEdges=size(ed2el,1); % Number of Edges
    Dbed= diag(n2ed(Dirichlet(:,1),Dirichlet(:,2)));  %Dirichlet boundary edges
    Nbed= diag(n2ed(Neumann(:,1),Neumann(:,2))); %Neumann boundary edges
    Redges=setdiff(1:nrEdges,Nbed); % removing Neumann edges from all edges
    L=length(Redges);  
    nE=size(elements,1);  % number of elements
    fprintf('number of elements = %d\n',nE);
    ndf=3*nE;            %degree of freedom       
    tic
    [B,D,M,b]=StiffMassConv_PH(coordinates,elements,params);
    fprintf('run-time(in sec) for the function StiffMassConv_PH.m  is = %d\n',toc);
    tic
    [C,e11,v,R11]= Lambda_PH(coordinates,ed2el,intE,Dbed,Redges,led,ledTN);
    fprintf('run-time(in sec) for the function Lambda_PH.m  is = %d\n',toc);
    %% Vectorization of Neumann BC
    LN=sparse(ndf,1); |\label{line:NeumannStart}|
    if (~isempty(Neumann))        
        nNb=size(Nbed,1); % number of Neumann boundary edges
        NP1=coordinates(Neumann(:,1),:); NP2=coordinates(Neumann(:,2),:); 
        Nmd=1/2.*(NP1+NP2); % mid-point of Neumann edges
        Nnorm=[NP2(:,2)-NP1(:,2) NP1(:,1)-NP2(:,1)]./vecnorm((NP1-NP2)')'; %normal to Neumann edge
        JN=reshape(((e11(Nbed)-ones(nNb,1))*3+[1 2 3])',3*nNb,1); 
        S=reshape(((diag(v(Nbed).*g(Nmd,nE,params)*Nnorm')).*R11(Nbed,:))',3*nNb,1);
        LN= accumarray (JN, S ,[ndf 1]); |\label{line:LN}|
    end     |\label{line:Neumannend}|
    %% Dirichlet boundary condition
    uhb1=sparse(L,1); |\label{line:DirichletStart}|
    if (~isempty(Dirichlet))
        [ed,~]=find(Redges==Dbed);  % extract the linear indices of Dirichlet edge number 
        me=(coordinates(Dirichlet(:,2),:)+coordinates(Dirichlet(:,1),:))./2; % mid-point of Dirichlet edges
        uhb1(ed)=-v(Dbed).*uD(me);
    end |\label{line:Dirichletend}|
%%%%%%% Solving saddle point problem %%%%%% 
    Mvec=[ B + D + M         -C'     |\label{line:AssemblyStart}|
         -C          sparse(L,L)];
    Fvec=[b+LN ; uhb1];   |\label{line:Assemblyend}|
    Wvec=Mvec\Fvec; 
    uh=Wvec(1:ndf);
    Kh=Wvec(1+ndf:length(Wvec));
    u=ue(coordinates(elements',:),nE,1); % exact solution
    K=exactLambda(elements,coordinates,ed2el,e11,led,Redges,params); |\label{line:exactLambda}|
    tic 
    [H1e(p1),L2e(p1),Merr(p1)]=PHerror(uh,C,B,Kh,K,grad4e,elements,coordinates); |\label{line:PHerror}|
    fprintf('run-time(in sec) for the function PHerror.m is = %d\n',toc);
end
%% order of convergence 
 for j=1:p1-1
    ocL2(j)=log(L2e(j)/L2e(j+1))/log(h(j)/h(j+1));
    ocH1(j)=log(H1e(j)/H1e(j+1))/log(h(j)/h(j+1));
    ocM(j)=log(Merr(j)/Merr(j+1))/log(h(j)/h(j+1));
 end
fprintf('order of convergences w.r.t h \n in   H1-norm   L2-norm   M-norm \n');
disp([ocH1',ocL2',ocM']);


\end{lstlisting}
\begin{itemize}

\item The computation of  $LN \in \mathbb{R}^{N \times 1}$ in MATLAB reads lines \ref{line:NeumannStart}-\ref{line:Neumannend}, incorporating the Neumann boundary conditions. The Neumann boundary condition in \eqref{eq:right g} at every mid-point $m_{E} \in \mathcal{N}_m(\Gamma_N)$ is evaluated.
\texttt{JN} and \texttt{S} $\in \mathbb{R}^{3 N_{n} \times 1}$ (\figurename{~\ref{fig:NeumannPH_diagram}}) indicate the row indexing vector and entries for the matrix $LN$.
Let $\eta_1, \eta_2, \dots, \eta_{N_{n}}$ be the Neumann boundary edge global numbers for the mesh $\mathcal{T}_h$.
\begin{figure}[H]
\centering
\includegraphics[width=11cm]{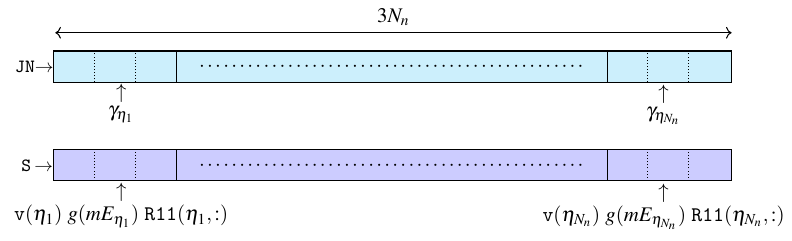}
\centering 
\caption{Geometrical representation of \texttt{JN}, and \texttt{S} vectors in MATLAB.}
\label{fig:NeumannPH_diagram}
\end{figure}
 The assembly of the matrix $LN$ is done in line \ref{line:LN} by using the $\mathbf{accumarray}$ with \texttt{JN} and \texttt{S}.


\item Line \ref{line:DirichletStart}-\ref{line:Dirichletend}: In taking care the Dirichlet boundary condition \eqref{eq:DirichletEqn}, the integral is approximated by the mid-point rule.

\item Line \ref{line:AssemblyStart}-\ref{line:Assemblyend}: Assembly of $\mathrm{M}_{vec}$ (see, \eqref{globalSystem}) in which zero ($L \times L$) matrix defined as a $\mathbf{sparse}$ matrix,  which helps to save run-time and storage space for zero-valued entries.

\item Line \ref{line:exactLambda}: The $\kappa$ \eqref{eq:exactLagrange} is calculated using the $\mathbf{exactLambda.m}$ function.

\item Line \ref{line:PHerror}: The errors  $u - u_h$ and $\kappa-\kappa_h$ over $X$ and $M$ are computed using the $\mathbf{PHerror.m}$ function.
\end{itemize} 
\begin{figure}[h]
\centering
\vspace{-0.4cm}
\includegraphics[width=1.05\textwidth]{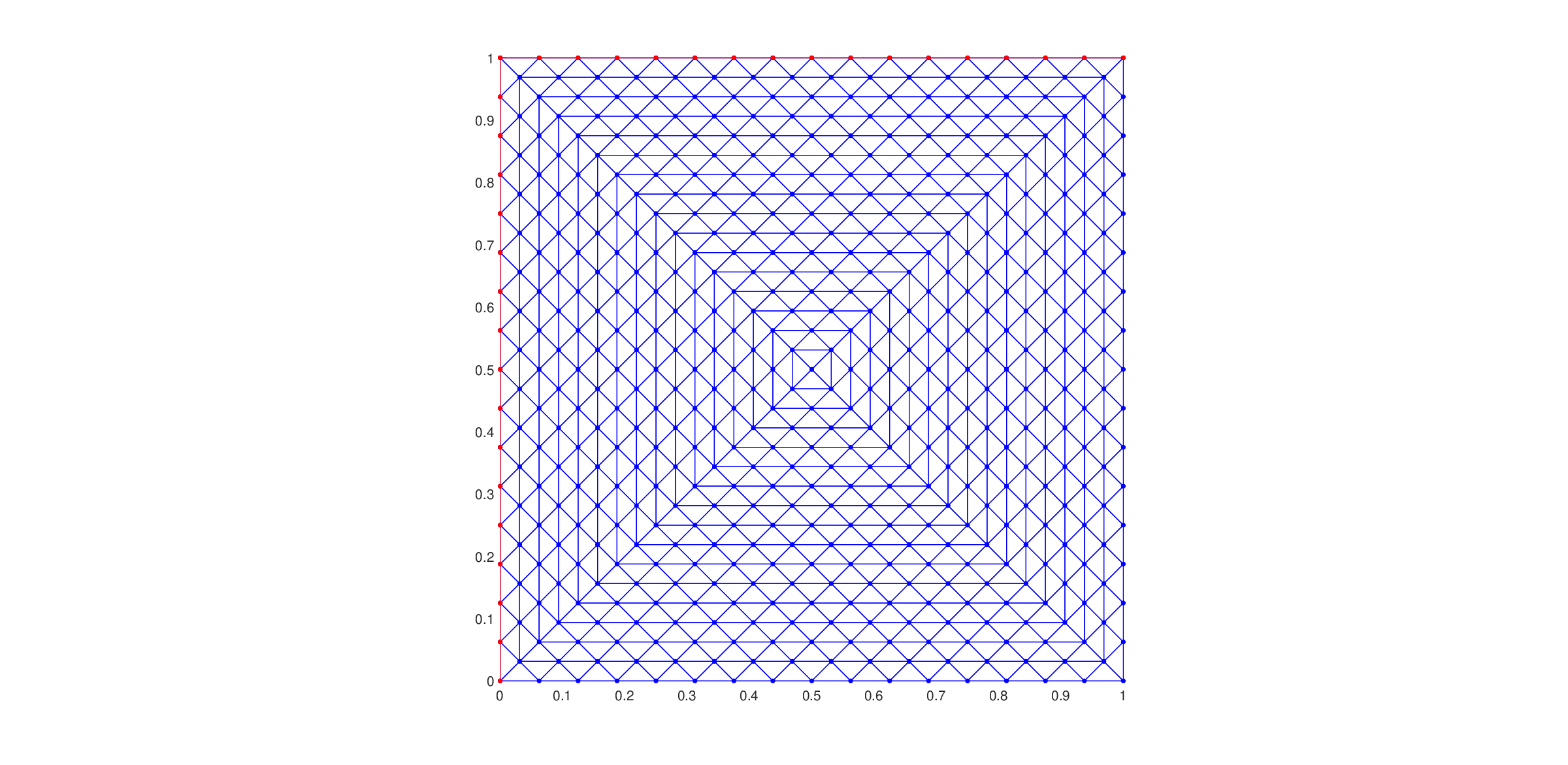} \\
\vspace{-0.4cm}
\caption{Example of (level 4, $h=1/16$) red-refined mesh used in computations of both elliptic and parabolic models. Red-colored sides denote $\Gamma_N$.}
\label{RedrefinedLevel4}
\end{figure}
\begin{example}
Let $~A=I_{2\times 2},~p=[1,1],~\delta=1 $. In \eqref{eq1}-\eqref{eq3} we take the exact solution  
$$u(x_1,x_2)=(x_1-x_1^2)(x_2-x_2^2)$$
defined on $\Omega = (0,1) \times (0,1)$ and corresponding to the load function $f$ in the form
\begin{equation*}
f(x_1,x_2)=2(x_2-x_2^2)+2(x_1-x_1^2)+(1-2x_1)(x_2-x_2^2)+(x_1-x_1^2)(1-2x_2)+u(x_1,x_2).
\end{equation*}
The Dirichlet boundary condition reads 
$$u_D=0 \quad \text{ on } \Gamma_D$$
and the Neumann boundary condition.
\begin{equation*}
g(x_1,x_2)=[(1-x_1-x_1^2)(x_2-x_2^2),~(x_1-x_1^2)(1-x_2-x_2^2)]\cdot \nu \quad \text{ on } \Gamma_N,
\end{equation*}
where $\nu$ denotes the outer normal vector (cf. Fig. \ref{RedrefinedLevel4}). \\

\begin{table}[ht]
\footnotesize
\captionsetup{font=footnotesize}
 \vspace{-0.1cm}
\begin{tabular*}{\textwidth}{@{\extracolsep\fill}crll}
 \vspace{-0.1cm}
\thead{Level} & \thead{\texttt{nE}} & \multicolumn{2}{@{}c@{}}{\thead{\textnormal{run-time (in sec)} }} \\   \cmidrule{3-4} %
 & & $\mathbf{StiffMassConv\_PH.m}$ & $\mathbf{Lambda\_PH.m}$ \\
\midrule
4 & 1024& 0.0029 &  0.0035  \\ 
 5 &  4096 & 0.01 & 0.017 \\   
  6 & 16384 & 0.018 & 0.025  \\   
7   & 65536 &  0.071 & 0.11  \\  
8  & 262144 & 0.3 & 0.51  \\  
9 & 1048576 & 1.23  & 2.13 \\
10 & 4194304 & 5.60 & 9.25\\
\end{tabular*}
\caption{Run-time of the $\mathbf{StiffMassConv\_PH.m}$ and $\mathbf{Lambda\_PH.m}$ for second order elliptic model problem.}
   \label{tab:StiffMassConv_Lambda_matrix_table}    
\end{table}

Table \ref{tab:StiffMassConv_Lambda_matrix_table} displays the runtime (in seconds) for vectorized implementations of functions $\mathbf{StiffMassConv\_PH.m}$ and $\mathbf{Lambda\_PH.m}$, which vary with the number of elements in the mesh. 

In \figurename{~\ref{linearRegression}}, we performed a linear regression analysis for the runtime of $\mathbf{StiffMassConv\_PH.m}$ (left) and $\mathbf{Lambda\_PH.m}$ (right) with respect to the number of elements. We obtained an R-squared approximately equal to 1 in both cases, concluding that our MATLAB implementations are of nearly linear time-scaling.

In Table \ref{tab:order_elliptic_table}, the errors $u-u_h$ in the $L^2$ and $H^1$-norm, and the error $\kappa-\kappa_h$ in $M$-norm over space M are presented.  In addition, we present the order of convergence in the respective norms.

\begin{table}[h]
\footnotesize
\captionsetup{font=footnotesize}
\begin{tabular*}{\textwidth}{@{\extracolsep\fill}rcccccc}
\thead{Level}  & \thead{$||u-u_h||_{H^1}$} & \thead{Order of \\ convergence} & \thead{$||u-u_h||_{L^2}$} & \thead{Order of \\ convergence} & \thead{$||\kappa-\kappa_h||_{M}$} & \thead{Order of \\ convergence}  \\  
\midrule
      1  & 0.0801 &    &  $0.0109$&    & 0.0232 &\\ 
      2  & 0.0397 & 1.0125   & $0.0025$ &   2.1050  &0.0107 &  1.1219 \\  
      3  & 0.0197  & 1.0083   & $6.22$e-04 & 2.0309   & 0.0055&  0.9678 \\  
      4  & 0.0099  & 1.0023    &  $1.54$e-04 & 2.0079  & 0.0028&  0.9866\\   
      5 & 0.0049  & 1.0006 & $3.86$e-05 &  2.0020   & 0.0014 & 0.9963 \\ 
    \end{tabular*}
    \caption{Order of convergence in $L^2$, $H^1$  and $M$-norms for second order elliptic model problem.}
   \label{tab:order_elliptic_table}
\end{table}
\begin{figure}[h!]
\small
\centering
\includegraphics[width=0.69\textwidth]{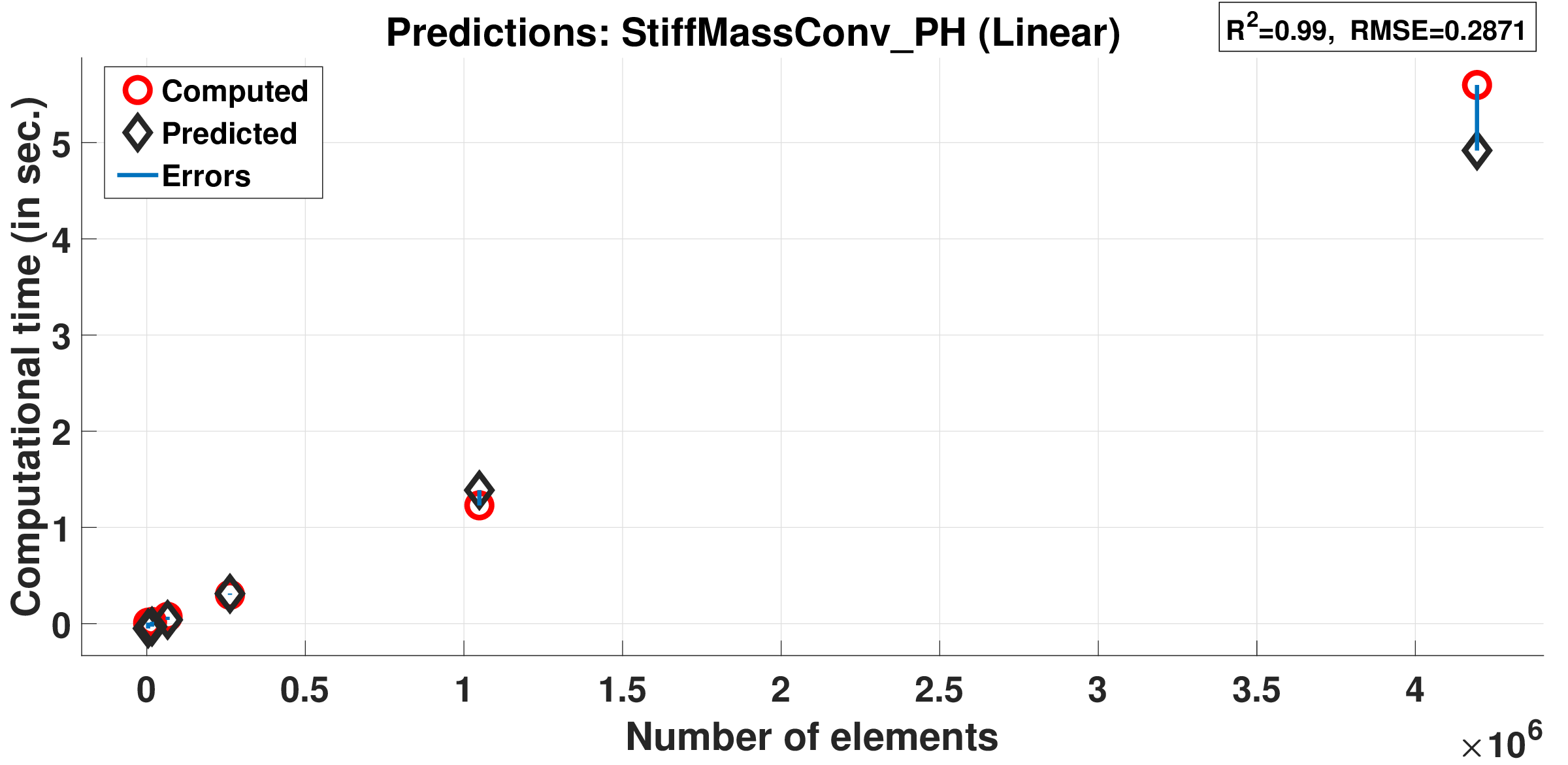} 
\includegraphics[width=0.69\textwidth]{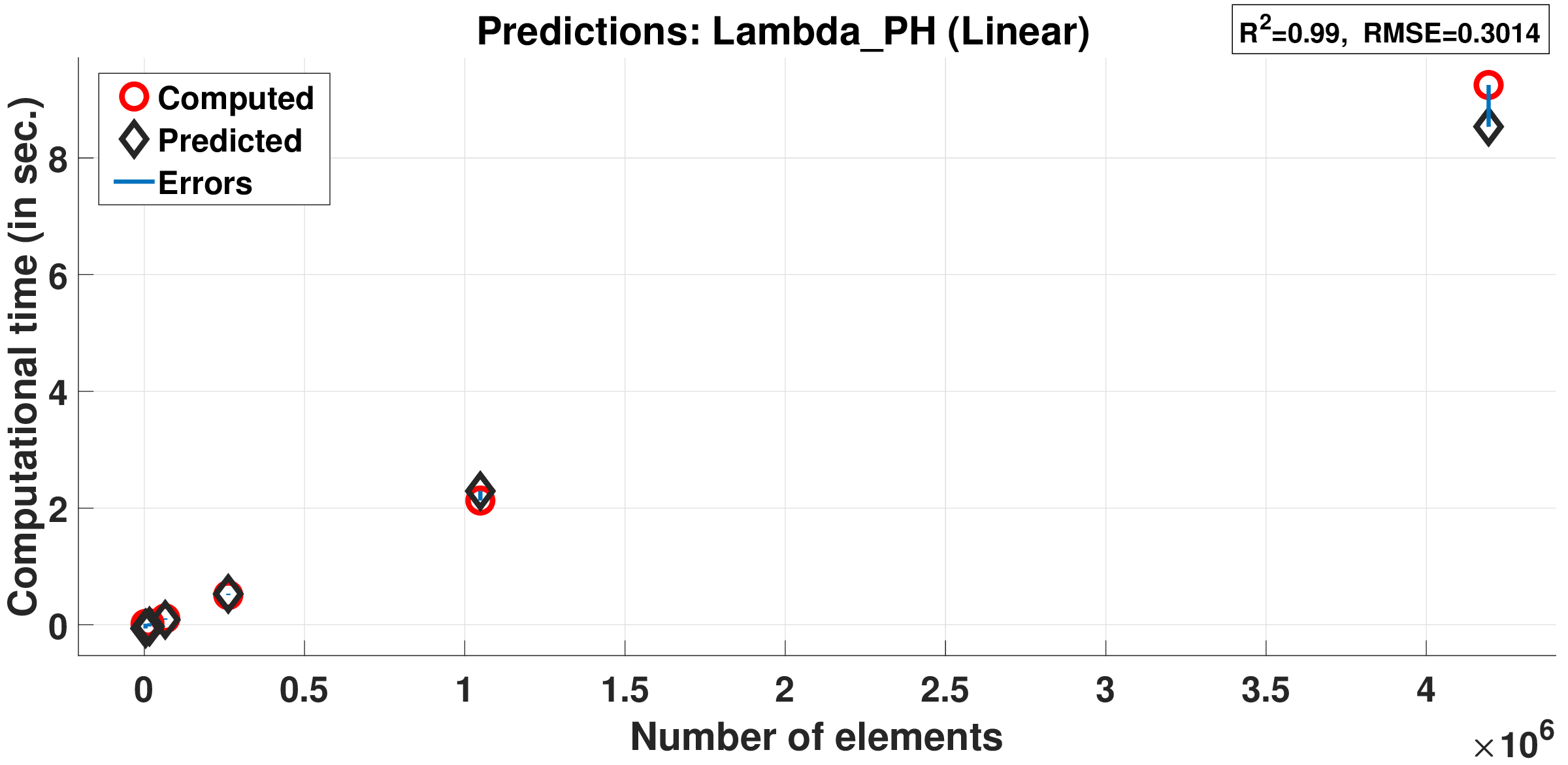}
\centering
\caption{Linear regression analysis between the computed and predicted values by Regression Learner (Statistics and Machine Learning Toolbox 11.7) in MATLAB.}
\label{linearRegression}
\end{figure}
\end{example}

\section{Implementation of  second order parabolic model problem}\label{sec:parabolicproblem}
\subsection{Parabolic problem and algebraic formulation
} Consider a second-order parabolic model problem with mixed boundary condition
\begin{align}
u_t(x,t)-\nabla\cdot(A\nabla u+up)+\delta u&=f(x,t)\hspace{0.3cm}\text{in}~\Omega\times (0,t_{\tilde N}],\label{paraeq1}\\
u(x,t)&=u_D(x,t)\hspace{0.3cm}\text{on}~\Gamma_D\times (0,t_{\tilde N}], \label{paraeq2} \\
 (A\nabla u+ub)\cdot\nu&=g(x,t)~~\text{on}~~\Gamma_N\times(0,t_{\tilde N}]\label{paraeq3}\\
u(x,0)&=u_{0}(x)\hspace{0.3cm}\text{in}~\Omega,\label{paraeq4}
\end{align}
where   $f$, $g$ and $u_{0}$ are appropriate smooth functions.
Primal hybrid formulation of  \eqref{paraeq1}-\eqref{paraeq4} is to seek a pair of solutions
$(u,\kappa):[0,t_{\tilde N}]\rightarrow X\times M$ such that
\begin{align}
(u_{t},v)+&a(u,v)+b(v,\kappa)=(f,v)+\langle g,v\rangle_{\Gamma_N},\hspace{0.2cm}\forall v\in X,\label{paraeq17}\\
&b(u,\chi)\qquad \qquad=\langle u_D,\chi\rangle,\hspace{0.7cm}\forall \chi\in M,\label{paraeq18}\\
&\qquad \qquad \quad u(0)=u_{0},\label{paraeq19}
\end{align}

Let the time interval $[0,t_{\tilde N}]$ be partitioned into equally spaced points $$0=t_0<t_1<\dots<t_{\tilde N}$$ with time step parameter $k$ such that $\tilde N=t_{\tilde N}/k$ and $t_{n}=nk$ for $n \in \{ 1, 2, \dots, \tilde N  \}$.
 For $\varphi\in \mathcal{C}[0,t_{\tilde N}]$,
we set the backward difference operator $\partial$ as 
$$\partial \varphi^{n}=(\varphi^{n}-\varphi^{n-1})/k$$
and for $t_{n-1/2}=(n-1/2)k$ the midpoint value as $$\varphi^{n-1/2}=(\varphi^{n}+\varphi^{n-1})/2.$$
 By employing the Crank-Nicolson implicit scheme in the time direction, we seek a solution pair
 $(u_h^n,\kappa_h^n)\in X_h\times M_h$
 satisfying
\begin{align*}
\sum_{T\in\mathcal{T}_h}\int_{T}\partial u_h^{n}v_{h}\dx+\sum_{T\in\mathcal{T}_h}\int_{T}\nabla u_h^{n-1/2}\cdot\nabla v_h\dx +\sum_{T\in\mathcal{T}_h}\int_{T}&(u_h^{n-1/2}p)\cdot\nabla v_h\dx \\
+\sum_{T\in\mathcal{T}_h}\int_{T} \delta u_h^{n-1/2}v_{h}\dx -\sum_{E\in\mathcal{E}^h/\mathcal{E}^h_{N}}\int_{E}\kappa_h^{n-1/2}[\![v_h]\!]_E\dgamma &=F^{t_{n-1/2}}(v_h),\\
-\sum_{E\in\mathcal{E}^h/\mathcal{E}^h_{N}}\int_{E}\chi_h[\![u_h^{n-1/2}]\!]_E\dgamma&=- \sum_{E\in\mathcal{E}^h/\mathcal{E}^h_{N}}\int_{E}\chi_hu^{n-1/2}_{Dh}\dgamma,
\end{align*}
for any testing pair $(v_h,\chi_h)\in  X_h\times M_h$. Here, 
$$F^{t_{n-1/2}}(v_h)=\displaystyle\sum_{T\in\mathcal{T}_h}\int_Tf(t_{n-1/2})v_h\dx+\sum_{E\in\mathcal{E}^{h}_{N}}\int_Eg(t_{n-1/2})v_h \dgamma .$$

\begin{remark}For this implicit scheme, we obtain a second-order convergence for the primal variable in $L^2$-norm (cf. \cite{thomee}), one can construct a complete discrete setup using any other implicit or explicit scheme.\end{remark}

Using the matrices $\mathbb M, \mathbb B, \mathbb D$, and $\mathbb C$ constructed for the elliptic case, we write the above algebraic formulation in matrix-vector form,
 \begin{equation}\label{eq:MatrixformParabolic}
 \text{M}_+~ \text{W}^{n} = \text{M}_-~ \text{W}^{n-1} + \text{F}.
\end{equation}
Here,
\begin{align*}
&\text{M}_{\pm}= \left(
\begin{array}{cc}
      (1 \pm \frac{1}{2}k \delta)\mathbb{M} \pm \frac{1}{2}k\mathbb{B} \pm \frac{1}{2}k\mathbb{D} & \quad \mp \frac{1}{2}k \mathbb{C}^{\prime}  \\
      \mp \frac{1}{2} \mathbb{C} & 0  \\
    \end{array}
  \right)_{(N+L)\times (N+L)}, \\
& \text{W}^{n} = \left(
    \begin{array}{c}
      U^n   \\
      \Lambda^n  \\
    \end{array}
  \right)_{(N+L)\times 1} , 
  \quad \text{W}^{n-1}= \left(
    \begin{array}{c}
      U^{n-1}   \\
      \Lambda^{n-1}  \\
    \end{array}
  \right)_{(N+L)\times 1}, \\
  & \text{F}= \left(
    \begin{array}{c}
      k F^{n-1/2}   \\
      b_D^{n-1/2}  \\
    \end{array}
  \right)_{(N+L)\times 1}, 
  \quad F^{n-1/2}=(F_j^{n-1/2})_{N\times 1}.
   \end{align*} 
For global basis functions $\phi_j$, 
   \begin{align*}
  &F_j^{n-1/2}
  = b^{n-1/2}_j + LN^{n-1/2}_j,
  \end{align*}
  where
  \begin{align*}
  b^{n-1/2}_j&\approx \displaystyle \sum_{T\in\mathcal{T}_h} \int_Tf^{n-1/2} \phi_j\dx,
  \quad  LN^{n-1/2}_j \approx \displaystyle \sum_{E\in\mathcal{E}^{h}_{N}} \int_Eg^{n-1/2} \phi_j\dgamma.
  \end{align*}
Then,
$b^{n-1/2}_j
= \frac{1}{2} b^{n-1}_j + \frac{1}{2} b^{n}_j
\approx \displaystyle \sum_{T\in\mathcal{T}_h} \frac{1}{2}\int_T \left(f^{n-1}\phi_j+f^{n}\phi_j\right) \dx.$
Here, $b^{n}_j$ is obtained as follows
\begin{align}\label{eq:parabolicRight_b}
b^{n}_j&\approx  \displaystyle \sum_{T\in\mathcal{T}_h} \frac{|T|}{3}~\sum_{i=1}^{3} f(m_{E_{i}},t_n)\phi_j(m_{E_{i}}).
\end{align}
The Neumann boundary term $LN^{n-1/2}_j = \displaystyle(LN^{n-1}_j+LN^{n}_j)/2$, and $LN^{n}_j$ is obtained by
\begin{equation}\label{eq:parabolicRight_g}
LN^{n}_j
\approx \displaystyle \sum_{E\in\mathcal{E}^{h}_{N}} |E| g(m_{E},t_n) \phi_j(m_{E}),~~~~~~\text{where}~ m_{E}\in \mathcal{N}_m(\Gamma_N).
\end{equation}
The component related to the Dirichlet boundary $b_D^{n-1/2}$ is defined as follows 
  \begin{align*}
  b_D^{n-1/2}=(b_{D_l}^{n-1/2})_{L\times 1},
  \quad 
  b_{D_l}^{n-1/2}\approx-\displaystyle\int_{E_l}\psi_l u^{n-1/2}_{Dh}\dgamma. 
\end{align*} 
Similarly, we obtain $b_{D_l}^{n-1/2}=(b_{D_l}^{n-1}+b_{D_l}^{n})/2$, and $b_{D_l}^{n}$ is obtained by
\begin{align}\label{eq:parabolicRight_uD}
b_{D_l}^{n}\approx-\displaystyle\int_{E_l}\psi_l u^{n}_{Dh}\dgamma &\approx -  |E_l| \sigma_l u^{n}_{Dh}(m_{E_l}),
\quad \text{where}~m_{E_l}\in \mathcal{N}_m(\Gamma_D) \nonumber \\ &= -|E_l| u^{n}_{Dh}(m_{E_l}), \quad \text{since}~\sigma_l=1~\text{for}~l=1,\dots,L.
\end{align}

The functions $\mathbf{StiffMassConv\_PH.m}$ and $\mathbf{Lambda\_PH.m}$ for the parabolic model problem will be similar to those for the elliptic model problem, which will handle the assembly of the matrices $\mathbb{M},~\mathbb{B},~\mathbb{C},~\mathbb{D}$. But, for the Crank-Nicolson scheme, we need to compute the vector $F^{n-1/2}$ at $t_{n-1/2}$, it is achieved by evaluating at $t_n$ and $t_{n-1}$ time points then taking an average of them. This necessitates the utilization of another function $\mathbf{load\_t.m}$ within $\mathbf{ParabolicMain.m}$. The $\mathbf{load\_t.m}$ calculates the values of $b^{n}$ \eqref{eq:parabolicRight_b} and $LN^{n}$ \eqref{eq:parabolicRight_g} at different time points $t_{n}$.

\begin{figure}
\centering
level=4 (h=1/16), time step: 8 of 16 \\
\includegraphics[width=0.6\textwidth]{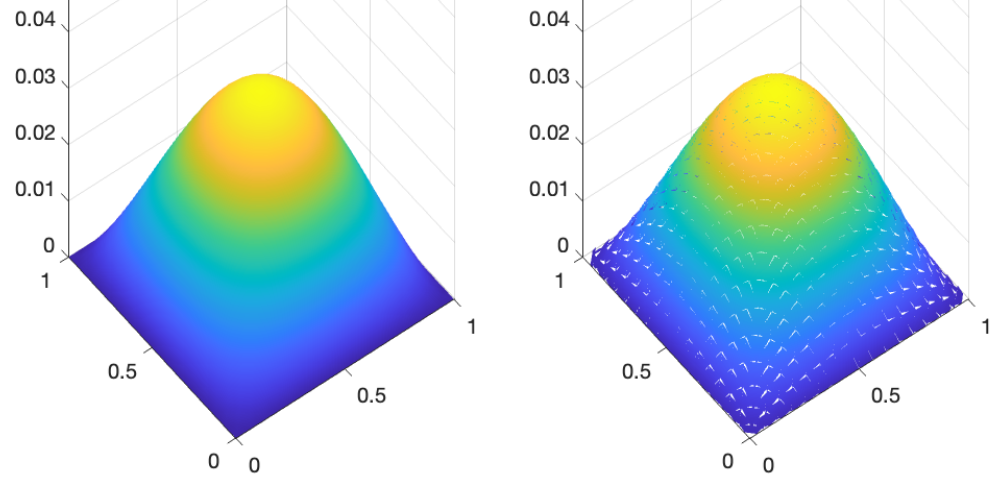} \\
\vspace{0.2cm}
level=4 (h=1/16), time step: 12 of 16 \\
\includegraphics[width=0.6\textwidth]{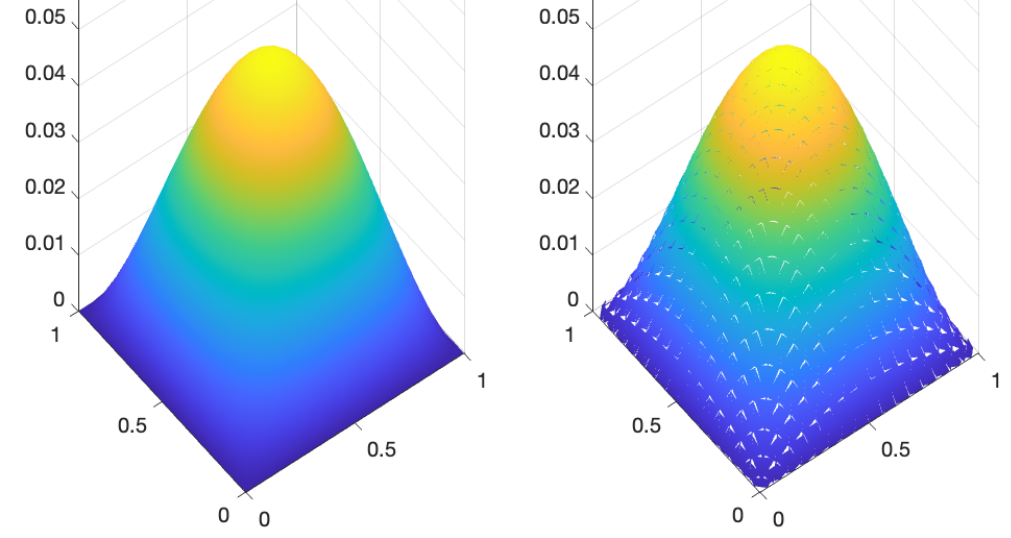} \\
\vspace{0.2cm}
level=4 (h=1/16), time step: 16 of 16 \\
\includegraphics[width=0.6\textwidth]{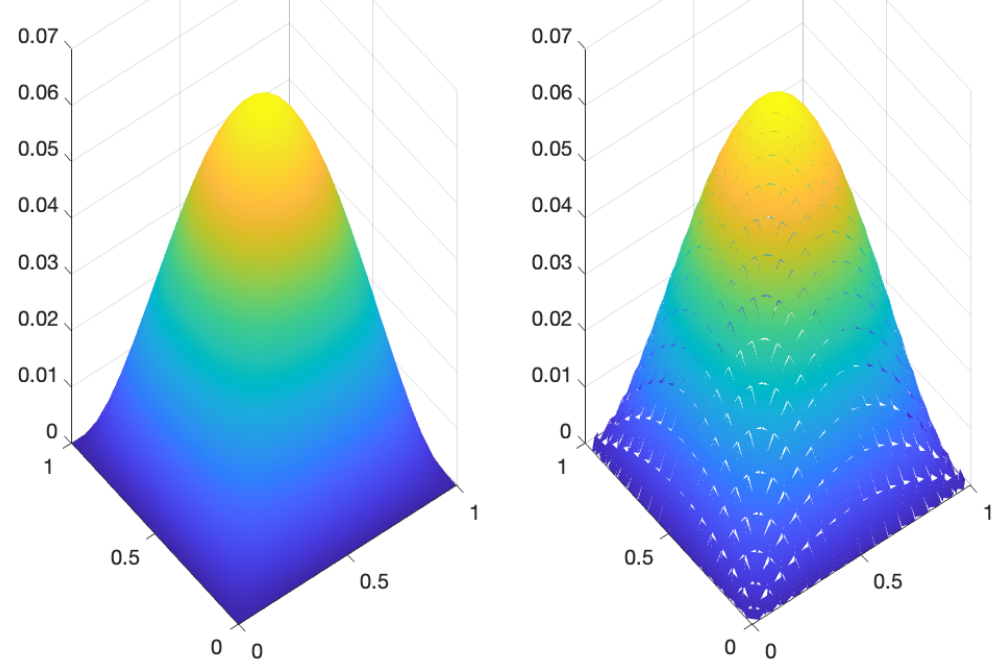} 
\caption{Exact $u$ (left column) vs. approximate $u_h$ solutions from the broken Sobolev space (right column) of the parabolic model at selected time steps.}\label{parabolic_u}
\centering
level=4 (h=1/16), time step: 16 of 16 \\
\includegraphics[width=0.6\textwidth]{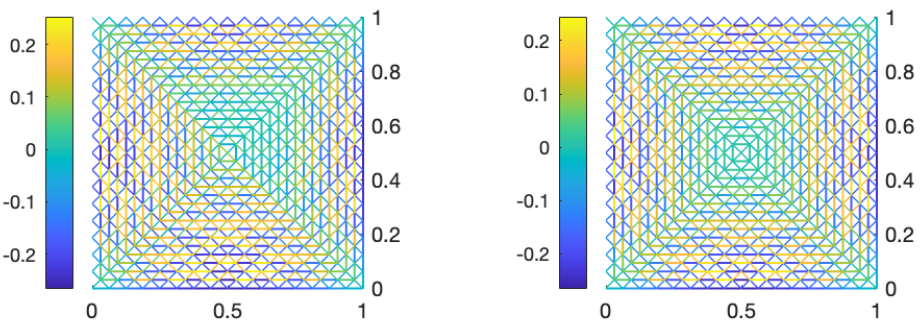} \\
\caption{Exact $\kappa$ (left column) and approximate $\kappa_h$ Lagrange multipliers (right column)  of the parabolic model in the last time step.}
\label{parabolic_lambda}
\end{figure}

\begin{example}

We consider again $A=I_{2\times 2},~ p=[1,1]~ \delta=1$ and the time interval $[0, 1]$.  In  \eqref{paraeq1}-\eqref{paraeq4} we take 
exact solution $u$ and corresponding load function $f$ in the forms
\begin{align*}
 u(x_1,x_2,t)=&t\,(x_1-x_1^2)(x_2-x_2^2),    \\
 f(x_1,x_2,t)=&(x_1-x_1^2)(x_2-x_2^2)+2t(x_2-x_2^2)+2t(x_1-x_1^2)\\
&+t(1-2x_1)(x_2-x_2^2)+t(x_1-x_1^2)(1-2x_2)+u(x_1,x_2,t).
\end{align*}
The Dirichlet and Neumann boundary (cf. Fig. \ref{RedrefinedLevel4} ) condition read $u_D=0$ 
and 
$$g(x_1,x_2)=[t(1-x_1-x_1^2)(x_2-x_2^2),~t(x_1-x_1^2)(1-x_2-x_2^2)]\cdot \nu~~
$$
Table \ref{tab:order_parabolic_table} shows the error $u-u_h$ in the $L^2$ and $H^1$-norms, and the error $\kappa-\kappa_h$ in the $M$-norm. For optimal order convergence of the Crank-Nicholson scheme, the space variable $h$ and the time parameter $k$ are considered equal, that is, $h=k$. Examples of particular exact and approximate solutions are shown in Figs. \ref{parabolic_u}, \ref{parabolic_lambda}.

\begin{table}[h]
\footnotesize
\captionsetup{font=footnotesize}
\begin{tabular*}{\textwidth}{@{\extracolsep\fill}rccccccc}
\thead{Level} & $h=k$ & \thead{$||u-u_h||_{H^1}$} & \thead{Order of \\ convergence} & \thead{$||u-u_h||_{L^2}$} & \thead{Order of \\ convergence} & \thead{$||\kappa-\kappa_h||_{M}$} & \thead{Order of \\ convergence}  \\  
\midrule
          1 & 1/2 & 0.1102 &    &  $0.0133$&    & 0.0397 &\\ 
      2 & 1/4 & 0.0493 & 1.1597   & $0.0028$ &   2.2400  &0.0131 &  1.5965 \\  
      3 & 1/8 & 0.0224  & 1.1366   & $6.70$e-04 & 2.0713 & 0.0061&  1.0977 \\  
      4 & 1/16 & 0.0106  & 1.0866   &  $1.59$e-04 & 2.0498  & 0.0029&  1.0702\\   
      5 & 1/32 & 0.0051 & 1.0483 & $3.88$e-05 &  2.0344   & 0.0014 & 1.0392 \\  
\end{tabular*}
    \caption{Order of convergence in $L^2$,$H^1$ and $M$-norms for second order parabolic model problem.}
   \label{tab:order_parabolic_table}
\end{table}
\end{example}

\section{Concluding remarks}\label{conclusion}
Primal hybrid FEM has been implemented efficiently in MATLAB for general two-dimensional elliptic and parabolic problems with mixed boundary conditions. Numerical experiments show that the codes are of nearly linear time-scaling, and the method converges optimally. The presented codes are flexible and may be extended to non-linear problems with an adequate amount of modification. We will include the implementations for the three-dimensional case in our future work.  

\section*{Acknowledgement}
 The second author's work is supported by
CSIR Extramural Research Grant. The third author announces the support of the Czech Science Foundation (GACR) through the GA23-04766S grant Variational approaches to dynamical problems in continuum mechanics.

\bibliography{main}


\end{document}